\documentclass{article}
\usepackage{amsmath}
\usepackage{amssymb}
\begin{document}
\def\e#1\e{\begin{equation}#1\end{equation}}
\def\ea#1\ea{\begin{align}#1\end{align}}
\def\eq#1{{\rm(\ref{#1})}}
\newtheorem{thm}{Theorem}[section]
\newtheorem{lem}[thm]{Lemma}
\newtheorem{prop}[thm]{Proposition}
\newtheorem{cor}[thm]{Corollary}
\newenvironment{dfn}{\medskip\refstepcounter{thm}
\noindent{\bf Definition \thesection.\arabic{thm}\ }}{\medskip}
\newenvironment{ex}{\medskip\refstepcounter{thm}
\noindent{\bf Example \thesection.\arabic{thm}\ }}{\medskip}
\newenvironment{proof}[1][,]{\medskip\ifcat,#1
\noindent{\it Proof.\ }\else\noindent{\it Proof of #1.\ }\fi}
{\hfill$\square$\medskip}
\def\dim{\mathop{\rm dim}}
\def\Re{\mathop{\rm Re}}
\def\Im{\mathop{\rm Im}}
\def\Image{\mathop{\rm Image}}
\def\Diff{\mathop{\rm Diff}}
\def\vol{\mathop{\rm vol}}
\def\Sym{{\ts\mathop{\rm Sym}}}
\def\GSym{{\ts\mathop{\rm GSym}}}
\def\Coad{\mathop{\rm Coad}}
\def\hcf{\mathop{\rm hcf}}
\def\Sp{\mathop{\rm Sp}}
\def\SO{\mathop{\rm SO}}
\def\U{\mathbin{\rm U}}
\def\SU{\mathop{\rm SU}}
\def\sn{{\ts\mathop{\rm sn}}}
\def\cn{{\ts\mathop{\rm cn}}}
\def\dn{{\ts\mathop{\rm dn}}}
\def\sech{\mathop{\rm sech}}
\def\cosech{\mathop{\rm cosech}}
\def\ge{\geqslant} 
\def\le{\leqslant} 
\def\cal{\mathcal}
\def\H{\mathbin{\mathbb H}}
\def\R{\mathbin{\mathbb R}}
\def\Z{\mathbin{\mathbb Z}}
\def\Q{\mathbin{\mathbb Q}}
\def\C{\mathbin{\mathbb C}}
\def\g{\mathbin{\mathfrak g}}
\def\CP{\mathbb{CP}}
\def\su{\mathfrak{su}}
\def\al{\alpha}
\def\be{\beta}
\def\ga{\gamma}
\def\de{\delta}
\def\ep{\epsilon}
\def\th{\theta}
\def\la{\lambda}
\def\vp{\varphi}
\def\La{\Lambda}
\def\Om{\Omega}
\def\Si{\Sigma}
\def\om{\omega}
\def\d{{\rm d}}
\def\pd{\partial}
\def\ts{\textstyle}
\def\w{\wedge}
\def\br{\buildrel}
\def\lt{\ltimes}
\def\sm{\setminus}
\def\ov{\overline}
\def\iy{\infty}
\def\ra{\rightarrow}
\def\t{\times}
\def\ha{{\ts{1\over2}}}
\def\op{\oplus}
\def\ti{\tilde}
\def\md#1{\vert#1\vert}
\def\ms#1{\vert#1\vert^2}
\def\bms#1{\bigl\vert#1\bigr\vert^2}
\def\bmd#1{\bigl\vert #1 \bigr\vert}
\def\an#1{\langle#1\rangle}
\def\ban#1{\bigl\langle#1\bigr\rangle}
\title{Special Lagrangian $m$-folds in $\C^m$ with symmetries}
\author{Dominic Joyce, \\ Lincoln College, Oxford}
\date{}
\maketitle

\section{Introduction}
\label{sy1}

This is the first in a series of papers constructing 
explicit examples of special Lagrangian submanifolds in $\C^m$. 
In it we will study special Lagrangian $m$-folds with large 
symmetry groups. The next five papers in the series are
\cite{Joyc2,Joyc3,Joyc4,Joyc5,Joyc6}. They use other methods 
to construct special Lagrangian $m$-folds, namely evolution 
equations, ruled submanifolds, and integrable systems.

The author's principal motivation for studying special Lagrangian
$m$-folds in $\C^m$ is that they provide local models for 
singularities of special Lagrangian $m$-folds in Calabi--Yau 
$m$-folds. In 1996 Strominger, Yau and Zaslow \cite{SYZ} proposed 
an explanation of mirror symmetry (the SYZ conjecture) between 
Calabi--Yau 3-folds $X,\hat X$ in terms of dual `fibrations' 
of $X$ and $\hat X$ by special Lagrangian $T^3$'s, with some 
singular fibres.

To make progress towards proving the SYZ conjecture, or even
stating it precisely, will require a good understanding of
the possible singularities that can develop in families of 
special Lagrangian 3-folds in a Calabi--Yau 3-fold. This
paper is part of a programme to develop such an understanding.

Some first steps in this direction were taken by the author 
in \cite{Joyc1}, which tried to define an invariant of 
Calabi--Yau 3-folds by counting special Lagrangian homology
3-spheres with weights; proving (or disproving) the conjectures 
made in \cite{Joyc1} will also require an understanding of the 
singularities of special Lagrangian 3-folds.

Perhaps the most obvious kinds of local model for singularities
of special Lagrangian $m$-folds are {\it special Lagrangian 
cones} in $\C^m$. The main results of the paper, in \S\ref{sy7} 
and \S\ref{sy8}, are a study of $\U(1)^{m-2}$-invariant special 
Lagrangian cones in $\C^m$, and the proof of the existence in 
$\C^m$ of large families of special Lagrangian cones on $T^{m-1}$ 
(for $m\ge 3$), ${\cal S}^2\t T^{m-3}$ (for $m\ge 4$) and 
${\cal S}^3\t T^{m-4}$ (for $m\ge 5$).

We begin in \S\ref{sy2} by introducing special Lagrangian
geometry in $\C^m$, and then in \S\ref{sy3} we give several
results relating to real analyticity of special Lagrangian
$m$-folds, including a construction of special Lagrangian
$m$-folds by evolving real analytic $(m\!-\!1)$-submanifolds
of $\C^m$. Section \ref{sy4} discusses {\it moment maps}, and 
shows that if $N$ is a Lagrangian submanifold with symmetry 
group $G$ then the moment map of $G$ is constant on~$N$. 

In \S\ref{sy5} we study {\it cohomogeneity one} SL $m$-folds $N$ in 
$\C^m$, where the orbits of the symmetry group $G\subset\SU(m)\lt\C^m$ 
are of codimension one in $N$. Then $N$ is foliated by a 1-parameter 
family of $G$-orbits parametrized by $t\in\R$. We write the condition 
that $N$ be special Lagrangian as an o.d.e.\ upon $G$-orbits depending 
on $t$, and by solving this equation we find examples of SL $m$-folds
in~$\C^m$.

Section \ref{sy6} considers {\it special Lagrangian cones} $N$ in 
$\C^m$. As cones are invariant under the group $\R_+$ of dilations 
of $\C^m$, by including dilations we can define the {\it generalized
symmetry group} $G\subset\R_+\t\SU(m)$ of $N$. So as in \S\ref{sy5} 
we can consider special Lagrangian cones on which the generalized 
symmetry group acts with cohomogeneity one.

As an example of this, in \S\ref{sy7} we study SL cones in $\C^m$
invariant under a subgroup $G\cong\U(1)^{m-2}$ of diagonal matrices
in $\C^m$, for $m\ge 3$. We reduce the problem to an o.d.e.\ in $m$
complex variables $w_1(t),\ldots,w_m(t)$, and by solving this o.d.e.\
fairly explicitly, we prove the existence of a large family of distinct
SL cones in $\C^m$ on $T^{m-1}$, and also of smaller families of SL
cones on ${\cal S}^2\t T^{m-3}$ in $\C^m$ for $m\ge 4$, and of SL 
cones on ${\cal S}^3\t T^{m-4}$ in $\C^m$ for~$m\ge 5$.

In \S\ref{sy8} we specialize to the case $m=3$, and consider 
$\U(1)$-invariant SL cones in $\C^3$ in more detail. Finally, section
\ref{sy9} gives a new construction of SL $m$-folds in $\C^m$ starting
with a special Lagrangian $m$-fold $L$ in $\C^m$ with `perpendicular 
symmetries', that is, vector fields in $\su(m)\lt\C^m$ which are
perpendicular to $L$ at every point.

We remark that Goldstein \cite[Th.~1]{Gold} has proved some 
results related to those of \S\ref{sy7}. He considers a compact
toric K\"ahler--Einstein $n$-fold $N$ with positive scalar 
curvature and a $\U(1)^n$-action preserving the structure. 
Then there is a unique, flat $\U(1)^n$-orbit $L$ which is 
minimal Lagrangian. Furthermore, he shows that there is at 
least one subgroup $\U(1)^{n-1}$ in $\U(1)^n$ with a sequence 
of non-flat minimal Lagrangian $\U(1)^{n-1}$-invariant tori 
$L_k$ converging to~$L$.

When Goldstein's results are applied to $\U(1)^{m-1}$ acting
on $\CP^{m-1}$ by isometries, they prove the existence of a 
family of $\U(1)^{m-2}$-invariant minimal Lagrangian tori 
$T^{m-1}$ in $\CP^{m-1}$, close to the unique minimal Lagrangian 
$\U(1)^{m-1}$-orbit. It can be shown that these lift to 
$\U(1)^{m-2}$-invariant special Lagrangian cones in $\C^m$, 
which is what we study in~\S\ref{sy7}.

There is also some overlap between the results of this paper, especially 
\S\ref{sy8}, and those of Castro and Urbano \cite{CaUr1,CaUr2} and 
Haskins \cite{Hask}. In particular, Haskins studies
$\U(1)$-invariant SL cones in $\C^3$, so that nearly all of \S\ref{sy8} 
is equivalent to results in \cite{Hask}, and also Theorem 
\ref{sy6thm2} below is essentially the same as 
\cite[Remark 1, p.~81--2]{CaUr2} and \cite[Th.~A]{Hask}. I would 
like to thank these authors for sending me copies of their work.
\bigskip

\noindent{\it Acknowledgements:} The author would like to thank
Nigel Hitchin, Mark Haskins, Karen Uhlenbeck, Ian McIntosh, 
Robert Bryant and Chuu-Lian Terng for helpful conversations,
and the referee for careful proofreading and suggesting
improvements.

\section{Special Lagrangian submanifolds in $\C^m$}
\label{sy2}

We begin by defining {\it calibrations} and {\it calibrated 
submanifolds}, following Harvey and Lawson~\cite{HaLa}. 

\begin{dfn} Let $(M,g)$ be a Riemannian manifold. An {\it oriented
tangent $k$-plane} $V$ on $M$ is a vector subspace $V$ of
some tangent space $T_xM$ to $M$ with $\dim V=k$, equipped
with an orientation. If $V$ is an oriented tangent $k$-plane
on $M$ then $g\vert_V$ is a Euclidean metric on $V$, so 
combining $g\vert_V$ with the orientation on $V$ gives a 
natural {\it volume form} $\vol_V$ on $V$, which is a 
$k$-form on~$V$.

Now let $\vp$ be a closed $k$-form on $M$. We say that
$\vp$ is a {\it calibration} on $M$ if for every oriented
$k$-plane $V$ on $M$ we have $\vp\vert_V\le \vol_V$. Here
$\vp\vert_V=\al\cdot\vol_V$ for some $\al\in\R$, and 
$\vp\vert_V\le\vol_V$ if $\al\le 1$. Let $N$ be an 
oriented submanifold of $M$ with dimension $k$. Then 
each tangent space $T_xN$ for $x\in N$ is an oriented
tangent $k$-plane. We say that $N$ is a {\it calibrated 
submanifold} or $\vp$-{\it submanifold} if 
$\vp\vert_{T_xN}=\vol_{T_xN}$ for all $x\in N$. 
\label{calibdef}
\end{dfn}

It is easy to show that calibrated submanifolds are automatically
{\it minimal submanifolds} \cite[Th.~II.4.2]{HaLa}. Here is 
the definition of special Lagrangian submanifolds in~$\C^m$.

\begin{dfn} Let $\C^m$ have complex coordinates $(z_1,\dots,z_m)$
and complex structure $I$, and define a metric $g$, a real 
2-form $\om$ and a complex $m$-form $\Om$ on $\C^m$ by
\begin{align*}
g=\ms{\d z_1}+\cdots+\ms{\d z_m},\quad
\om&={i\over 2}(\d z_1\w\d\bar z_1+\cdots+\d z_m\w\d\bar z_m),\\
\text{and}\quad\Om&=\d z_1\w\cdots\w\d z_m.
\end{align*}
Then $\Re\Om$ and $\Im\Om$ are real $m$-forms on $\C^m$. Let
$L$ be an oriented real submanifold of $\C^m$ of real dimension 
$m$, and let $\th\in[0,2\pi)$. We say that $L$ is a {\it special 
Lagrangian submanifold} of $\C^m$ with {\it phase} ${\rm e}^{i\th}$,
if $L$ is calibrated with respect to $\cos\th\,\Re\Om+\sin\th\,\Im\Om$,
in the sense of Definition~\ref{calibdef}.

We will often abbreviate `special Lagrangian' by `SL', and 
`$m$-dimensional submanifold' by `$m$-fold', so that we shall
talk about SL $m$-folds in $\C^m$. Usually we take $\th=0$, so 
that $L$ has phase 1, and is calibrated with respect to $\Re\Om$. 
When we discuss special Lagrangian submanifolds without specifying 
a phase, we mean them to have phase~1.
\label{cmsldef}
\end{dfn}

Harvey and Lawson \cite[Cor.~III.1.11]{HaLa} give the following 
alternative characterization of special Lagrangian submanifolds.

\begin{prop} Let\/ $L$ be a real $m$-dimensional submanifold 
of $\C^m$. Then $L$ admits an orientation making it into an
SL submanifold of\/ $\C^m$ with phase ${\rm e}^{i\th}$ if and 
only if\/ $\om\vert_L\equiv 0$ 
and\/~$(\sin\th\,\Re\Om-\cos\th\,\Im\Om)\vert_L\equiv 0$.
\label{sl1prop}
\end{prop}

Note that an $m$-dimensional submanifold $L$ in $\C^m$ is 
called {\it Lagrangian} if $\om\vert_L\equiv 0$. Thus special 
Lagrangian submanifolds are Lagrangian submanifolds satisfying 
the extra condition that $(\sin\th\,\Re\Om-\cos\th\,\Im\Om)\vert_L
\equiv 0$, which is how they get their name.

\section{Real analyticity of SL submanifolds} 
\label{sy3}

In this section we collect together several results about
special Lagrangian submanifolds in $\C^m$ related to {\it real
analyticity}. For simplicity we fix the phase of all special 
Lagrangian submanifolds to be 1. As special Lagrangian 
submanifolds in $\C^m$ are calibrated, they are locally 
minimal. Harvey and Lawson \cite[Th.~III.2.7]{HaLa} use 
this to show that they are real analytic:

\begin{thm} Let\/ $L$ be a special Lagrangian submanifold in
$\C^m$. Then $L$ is real analytic wherever it is nonsingular.
\label{morrthm}
\end{thm}

Note that the restriction to nonsingular $L$ is necessary here,
as Harvey and Lawson \cite[p.~97]{HaLa} give examples of
singularities of SL submanifolds in $\C^m$ which
are not real analytic. Harvey and Lawson \cite[Th.~III.5.5]{HaLa} 
also use real analyticity in a different way, to prove the following 
result.

\begin{thm} Suppose $P$ is a real analytic $(m-1)$-submanifold
of\/ $\C^m$ with\/ $\om\vert_P\equiv 0$. Then there exists a 
locally unique special Lagrangian submanifold\/ $N$ of\/ $\C^m$ 
containing~$P$.
\label{hlthm}
\end{thm}

They assume $P$ is real analytic because their proof uses the
{\it Cartan--K\"ahler Theorem}, from the subject of exterior 
differential systems, and this only works in the real analytic
category. One can think of the submanifold $N$ as defined by a 
kind of Taylor series, which converges in a small neighbourhood 
of~$P$. 

We will now show that $N$ is the total space of a 1-parameter 
family $\bigl\{P_t:t\in(-\ep,\ep)\bigr\}$ of real analytic 
submanifolds of $\C^m$ diffeomorphic to $P$, which satisfy a 
first-order o.d.e.\ in $t$, with initial data $P_0=P$. This 
will provide motivation for several of the constructions of 
special Lagrangian submanifolds in $\C^m$ to be considered 
in this paper and its sequels.

\begin{thm} Let\/ $P$ be a compact, orientable, real analytic 
$(m-1)$-manifold, $\chi$ a real analytic, nowhere vanishing 
section of\/ $\La^{m-1}TP$, and\/ $\phi:P\ra\C^m$ a real analytic
embedding (immersion) such that\/ $\phi^*(\om)\equiv 0$ on 
$P$. Then there exists $\ep>0$ and a unique family 
$\bigl\{\phi_t:t\in(-\ep,\ep)\bigr\}$ of real analytic maps 
$\phi_t:P\ra\C^m$ with $\phi_0=\phi$, satisfying the equation
\e
\left({\d\phi_t\over\d t}\right)^c=(\phi_t)_*(\chi)^{b_1\ldots b_{m-1}}
(\Re\Om)_{b_1\ldots b_{m-1}b_m}g^{b_mc},
\label{dphiteq}
\e
using the index notation for (real) tensors on $\C^m$, where
$g^{bc}$ is the inverse of the Euclidean metric on $\C^m$. Define 
$\Phi:(-\ep,\ep)\t P\ra\C^m$ by $\Phi(t,p)=\phi_t(p)$. Then 
$N=\Image\Phi$ is a nonsingular embedded (immersed) special 
Lagrangian submanifold of\/~$\C^m$.
\label{slevthm}
\end{thm}

\begin{proof} Equation \eq{dphiteq} is an evolution equation for 
the maps $\phi_t:P\ra\C^m$, with initial condition $\phi_0=\phi$.
As $P$ is compact and everything is real analytic, the existence 
of a unique solution for $t$ in $(-\ep,\ep)$ for some $\ep>0$
follows from standard techniques in partial differential equations.
For instance, one can prove it by applying the Cauchy--Kowalevsky 
Theorem \cite[p.~234]{Rack} to an evolution equation for
$(\phi_t,\d\phi_t)$ derived from~\eq{dphiteq}.

Thus the family $\bigl\{\phi_t:t\in(-\ep,\ep)\bigr\}$ exists, 
and it remains to prove that $N=\Image\Phi$ is special Lagrangian.
Now by Theorem \ref{hlthm}, as $\phi^*(\om)\equiv 0$ and
$\phi$ is real analytic, there is a locally unique real analytic 
special Lagrangian submanifold $N'$ in $\C^m$ containing 
$\phi(P)$. We shall show that $N'=N$ locally.

To do this, observe that \eq{dphiteq} also makes sense as
an evolution equation for submanifolds of $N'$. That is, we 
could look for a family $\bigl\{\phi_t':t\in(-\ep',\ep')\bigr\}$ 
of real analytic maps $\phi_t':P\ra N'$ with $\phi_0=\phi$, 
satisfying the equation
\e
\left({\d\phi_t'\over\d t}\right)^c=
(\phi_t')_*(\chi)^{b_1\ldots b_{m-1}}
(\Re\Om\vert_{N'})_{b_1\ldots b_{m-1}b_m}(g\vert_{N'})^{b_mc},
\label{sy3eq1}
\e
using the index notation for tensors on $N'$. It follows as above
that for some $\ep'>0$ there exists a unique solution to this 
problem. 

Let $p\in P$ and $t\in(-\ep',\ep')$, and set $x=\phi_t'(p)$, so
that $x\in N'$. Treating $\C^m$ and $(\C^m)^*$ as {\it real}\/ 
vector spaces, we have orthogonal direct sums $\C^m=T_xN'\op V$ 
and $(\C^m)^*=T_x^*N'\op V^*$, where $V$ is the perpendicular
subspace to $T_xN'$. This induces a splitting~$\La^m(\C^m)^*=
\bigoplus_{k=0}^m\La^kT_x^*N'\otimes\La^{m-k}V^*$.

Now $\Re\Om\in\La^m(\C^m)^*$, and $N'$ is calibrated with respect
to $\Re\Om$. This implies that the component of $\Re\Om$ in
$\La^{m-1}T_x^*N'\otimes V^*$ is zero, because this measures the 
change in $\Re\Om\vert_{T_xN'}$ under small variations of the 
subspace $T_xN'$, but $\Re\Om\vert_{T_xN'}$ is maximum and 
therefore stationary.

Since $(\phi_t')_*(\chi)\vert_p$ lies in $\La^{m-1}T_xN'$,
it follows that 
\begin{equation*}
(\phi_t')_*(\chi)^{b_1\ldots b_{m-1}}\big\vert_p
(\Re\Om)_{b_1\ldots b_{m-1}b_m}
\in T_x^*N'\subset(\C^m)^*,
\end{equation*}
as the component in $V^*$ comes from the component of $\Re\Om$ in 
$\La^{m-1}T_x^*N'\otimes V^*$, which is zero. Therefore
\e
\begin{split}
&(\phi_t')_*(\chi)^{b_1\ldots b_{m-1}}\big\vert_p
(\Re\Om)_{b_1\ldots b_{m-1}b_m}=\\
&(\phi_t')_*(\chi)^{b_1\ldots b_{m-1}}\big\vert_p
(\Re\Om\vert_{T_xN'})_{b_1\ldots b_{m-1}b_m}.
\end{split}
\label{sy3eq2}
\e

Because the splitting $(\C^m)^*=T_x^*N'\op V^*$ is orthogonal we have 
$g^{bc}=(g\vert_{T_xN'})^{bc}+h^{bc}$ for some $h\in S^2V$. Contracting
this with \eq{sy3eq2} shows that
\e
\begin{split}
&(\phi_t')_*(\chi)^{b_1\ldots b_{m-1}}
(\Re\Om)_{b_1\ldots b_{m-1}b_m}g^{b_mc}=\\
&(\phi_t')_*(\chi)^{b_1\ldots b_{m-1}}
(\Re\Om\vert_{N'})_{b_1\ldots b_{m-1}b_m}(g\vert_{N'})^{b_mc}
\end{split}
\label{sy3eq3}
\e
at $p$, as the r.h.s.\ of \eq{sy3eq2} lies in $T_x^*N'$, so its
contraction with $h$ is zero. 

Equation \eq{sy3eq3} holds for all $p\in P$ and $t\in(-\ep',\ep')$.
Therefore by \eq{sy3eq1} the $\phi_t'$ also satisfy \eq{dphiteq}, 
and so $\phi_t'=\phi_t$ by uniqueness. Hence $\phi_t$ maps $P$ to 
$N'$, and $\Phi$ maps $(-\ep,\ep)\t P$ to $N'$, if $\ep$ is 
sufficiently small.

Finally, suppose $\phi=\phi_0$ is an embedding. Then $\phi_t:P\ra N'$ 
is also an embedding for small $t$. But $\d\phi_t/\d t$ is a normal 
vector field to $\phi_t(P)$ in $N'$, with length $\bmd{(\phi_t)_*(\chi)}$.
As $\chi$ is nonvanishing, this vector field is nonzero, so $\Phi$ is 
an embedding for small $\ep$, with $\Image\Phi$ an open subset of $N'$, 
which is therefore special Lagrangian. If $\phi$ is an immersion,
then $\Phi$ is a special Lagrangian immersion, in a similar way.
\end{proof}

The condition that $P$ be compact is not always necessary here.
Whether $P$ is compact or not, in a small neighbourhood of any 
$p\in P$ the maps $\phi_t$ always exist for $t\in(-\ep,\ep)$ and 
some $\ep>0$, which may depend on $p$. If $P$ is compact we can 
choose an $\ep>0$ valid for all $p$, but if $P$ is noncompact
there may not exist such an~$\ep$.

We can also relax the condition that $\phi:P\ra\C^m$ be an 
embedding or an immersion, and instead require only that $\phi$ 
be real analytic. Then the conclusions of the theorem still hold, 
except that $\Phi$ is no longer an embedding or an immersion, and 
$\Image\Phi$ will in general be a {\it singular} special Lagrangian
submanifold of $\C^m$. This can be used as a technique for
constructing singular special Lagrangian submanifolds.

\section{Symmetries and moment maps}
\label{sy4}

Let $\C^m$ have its usual metric $g$ and K\"ahler form $\om$. Then 
the group of automorphisms of $\C^m$ preserving $g$ and $\om$ is 
$\U(m)\lt\C^m$, where $\C^m$ acts by translations. Let $G$ be a 
Lie subgroup of $\U(m)\lt\C^m$, with Lie algebra $\g$, and let 
$\phi:\g\ra C^\iy(T\C^m)$ be the natural action of $\g$ on $\C^m$ 
by vector fields. 

Then a {\it moment map} for the action of $G$ on $\C^m$ is a 
smooth map $\mu:\C^m\ra\g^*$, such that 
\begin{itemize}
\item[(a)] $\iota(\phi(x))\om=x\cdot\d\mu$ for all $x\in\g$, where
`$\cdot$' is the pairing between $\g$ and~$\g^*$,
\item[(b)] $\mu$ is equivariant with respect to the $G$-action
$\Phi$ on $\C^m$ and the coadjoint $G$-action on $\g^*$.
\end{itemize}
If $G$ is compact or semisimple then a moment map $\mu$ always exists,
but in general there may be obstructions to the existence of~$\mu$.

The subsets $\mu^{-1}(c)$ for $c\in\g^*$ are called {\it level sets} 
of the moment map. Define the {\it centre} $Z(\g^*)$ to be the
vector subspace of $\g^*$ fixed by the coadjoint action of $G$.
Then, as $\mu(\ga\cdot z)=\Coad(\ga)\mu(z)$ for each $z\in M$ and
$\ga\in G$, we see that $\mu^{-1}(c)$ is $G$-invariant if and only
if~$c\in Z(\g^*)$.

Here is a result characterizing $G$-orbits $\cal O$ 
with~$\om\vert_{\cal O}\equiv 0$.

\begin{prop} Let\/ $G$ be a connected Lie subgroup of\/ $\U(m)\lt\C^m$ 
with Lie algebra $\g$ and moment map $\mu:\C^m\ra\g^*$, and let\/ $\cal O$ 
be an orbit of\/ $G$ in $\C^m$. Then $\om\vert_{\cal O}\equiv 0$ if and 
only if\/ ${\cal O}\subseteq\mu^{-1}(c)$ for some~$c\in Z(\g^*)$.
\label{sy4prop1}
\end{prop}

\begin{proof} Let $x,y\in\g$, and let $\phi(x),\phi(y)$ be the
induced vector fields on $\C^m$. Then, by definition of the
moment map $\mu$, we see that
\e
\om\bigl(\phi(x),\phi(y)\bigr)=\phi(y)\cdot(x\cdot\d\mu)=
x\cdot\bigl(\phi(y)\cdot\d\mu\bigr)=x\cdot{\cal L}_{\phi(y)}\mu,
\label{om0xyeq}
\e
where ${\cal L}_{\phi(y)}$ is the Lie derivative. As $\cal O$ is a 
$G$-orbit, the vector fields $\phi(x)$ for $x\in\g$ are tangent to 
$\cal O$, and generate its whole tangent space. Thus, $\om\vert_{\cal O}
\equiv 0$ if and only if $\om\bigl(\phi(x),\phi(y)\bigr)=0$ on 
$\cal O$ for all~$x,y\in\g$. 

By \eq{om0xyeq} this holds if and only if ${\cal L}_{\phi(y)}\mu=0$ 
on $\cal O$ for all $y\in\g$. Since $G$ is connected, this is true 
if and only if $\mu$ is constant on $\cal O$, so that 
${\cal O}\subseteq\mu^{-1}(c)$ for some $c\in\g^*$. But $\cal O$ is 
$G$-invariant, so $c$ is $\Coad(G)$-invariant, and~$c\in Z(\g^*)$.
\end{proof}

The reason we are interested in moment maps is that $G$-invariant 
Lagrangian submanifolds in $\C^m$ lie in level sets of the moment map 
$\mu$ of $G$. Thus, moment maps are a tool for studying Lagrangian 
(and hence special Lagrangian) submanifolds with symmetries.

\begin{prop} Suppose $N$ is a connected Lagrangian submanifold 
of\/ $\C^m$, and\/ $G$ a connected Lie subgroup of\/ $\U(m)\lt\C^m$ 
preserving $N$, with Lie algebra $\g$. Then $G$ admits a moment
map $\mu$, and\/ $N\subseteq\mu^{-1}(c)$ for some~$c\in Z(\g^*)$.
\label{sy4prop2}
\end{prop}

\begin{proof} Let $x\in\g$. Then $\iota(\phi(x))\om$ is a closed 
1-form on $\C^m$, so there exists a smooth function $f_x:\C^m\ra\R$, 
unique up to addition of a constant, with $\d f_x=\iota(\phi(x))\om$. 
Since $\om\vert_N\equiv 0$ and $\phi(x)$ is tangent to $N$ we have 
$\d f_x\vert_N\equiv 0$, so $f_x$ is constant on $N$, as $N$ is 
connected.

As the functions $f_x$ are defined up to addition of a
constant, we can choose $f_x$ uniquely by requiring that
$f_x=0$ on $N$. Clearly $f_x$ is linear in $x$. Hence there 
is a unique map $\mu_0:\C^m\ra\g^*$ with $x\cdot\mu_0=f_x$ for 
all $x\in\g$, and thus $\iota(\phi(x))\om=x\cdot\d\mu_0$ as in 
part (a) above. The $G$-equivariance in part (b) follows because 
$N$ is $G$-invariant.

Thus we have constructed a particular moment map $\mu_0$ for $G$, 
with the property that $\mu_0\equiv 0$ on $N$. Hence $G$ admits a 
moment map. If $\mu$ is any moment map for $G$ then $\mu-\mu_0$ is 
a constant $c\in Z(\g^*)$, and so~$N\subseteq\mu^{-1}(c)$.
\end{proof}

Now consider special Lagrangian $m$-folds $N$ in $\C^m$. Then
instead of $\U(m)\lt\C^m$ we should use $\SU(m)\lt\C^m$, the group 
of automorphisms of $\C^m$ preserving $g,\om$ and~$\Om$.

\begin{dfn} Let $N$ be a special Lagrangian $m$-fold in $\C^m$.
Define the {\it symmetry group} $\Sym(N)$ of $N$ to be the 
Lie subgroup of $\SU(m)\lt\C^m$ preserving $N$. That is, the 
elements of $\Sym(N)$ are automorphisms of $\C^m$ preserving $g$, 
$\om$, $\Om$ and $N$. Define the {\it restricted symmetry group} 
$\Sym^0(N)$ of $N$ to be the connected component of $\Sym(N)$
containing the identity. 
\end{dfn}

From Proposition \ref{sy4prop2} we get:

\begin{cor} Let\/ $N$ be a connected special Lagrangian $m$-fold 
in $\C^m$, and set\/ $G=\Sym^0(N)\subset\SU(m)\lt\C^m$. Then $G$ 
admits a moment map $\mu:\C^m\ra\g^*$, and\/ $N\subseteq\mu^{-1}(c)$ 
for some~$c\in Z(\g^*)$.
\end{cor}

\subsection{Special Lagrangian $m$-folds with cohomogeneity one}
\label{sy41}

Let $N$ be a special Lagrangian submanifold of $\C^m$. The 
symmetry group $\Sym(N)$ was defined above to be the Lie 
subgroup of $\SU(m)\ltimes\C^m$ preserving $N$. Now it is a 
general principle that the easiest geometric objects to 
construct are those with large symmetry groups. It can be shown
that all {\it homogeneous} special Lagrangian submanifolds are 
affine subspaces $\R^m$ in $\C^m$, which are not very interesting.

The next most symmetric kinds of special Lagrangian submanifold $N$ 
are those of {\it cohomogeneity one}, that is, where the orbits of 
the symmetry group are of codimension one in $N$. For these we can
prove the following theorem.

\begin{thm} Let\/ $G$ be a connected Lie subgroup of\/ 
$\SU(m)\ltimes\C^m$ with Lie algebra $\g$ and moment map 
$\mu:\C^m\ra\g^*$, let\/ $\cal O$ be an oriented orbit of\/ 
$G$ in $\C^m$ with\/ $\dim{\cal O}=m-1$, and suppose 
${\cal O}\subset\mu^{-1}(c)$ for some $c\in Z(\g^*)$. Then 
there exists a locally unique, $G$-invariant special Lagrangian 
submanifold\/ $N$ in $\C^m$ containing $\cal O$. Furthermore 
$N\subset\mu^{-1}(c)$, and\/ $N$ is fibred by $G$-orbits 
isomorphic to $\cal O$ near $\cal O$. Thus $N$ is locally 
diffeomorphic to $(-\ep,\ep)\t{\cal O}$, for some $\ep>0$, 
and we can think of\/ $N$ as a smooth curve of\/ $G$-orbits.
\label{sy4thm}
\end{thm}

\begin{proof} Proposition \ref{sy4prop1} shows that 
$\om\vert_{\cal O}\equiv 0$, and as $\cal O$ is a $G$-orbit
it is a real analytic submanifold of $\C^m$, with dimension
$m-1$ by assumption. Therefore Theorem \ref{hlthm} shows that
there exists a locally unique SL submanifold
$N$ in $\C^m$ containing $\cal O$. But $\cal O$ is $G$-invariant 
and $G\subset\SU(m)\ltimes\C^m$, so $N$ must also be (locally) 
$G$-invariant, by local uniqueness. Hence $N\subseteq\mu^{-1}(c)$, 
by Proposition~\ref{sy4prop2}.

We can also construct $N$ by an evolution equation, as in Theorem
\ref{slevthm}. Choose $x\in{\cal O}$, and let $H$ be the stabilizer
of $x$ in $G$. Set $P=G/H$, and define $\phi:P\ra\C^m$ by
$\phi(\ga H)=\ga x$ for all $\ga\in G$. Then $\phi$ is an
immersion, with $\phi(P)={\cal O}$. Clearly $P$ and $\phi$ are 
real analytic, and $\phi^*(\om)\equiv 0$ on~$P$.

As $P\cong{\cal O}$ is oriented and $G$ acts on it by isometries, 
we can choose a nonvanishing, $G$-invariant section $\chi$ of
$\La^{m-1}TP$. Suppose for the moment that $P$ is compact.
Then Theorem \ref{slevthm} applies to give an embedding 
$\Phi:(-\ep,\ep)\t P\ra\C^m$, whose image is an open subset of 
$N$. As $\chi$ is $G$-invariant and $\phi$ is $G$-equivariant,
we see by uniqueness that $\Phi$ is equivariant under the actions 
of $G$ on $P$ and~$\C^m$. 

Thus $\Phi\bigl(\{t\}\t P\bigr)$ is a $G$-orbit in $\C^m$ isomorphic 
to $\cal O$ for $t\in(-\ep,\ep)$, and $N$ is fibred by a smooth 
1-parameter family of $G$-orbits isomorphic to $\cal O$ near $\cal O$.
This completes the proof, except that we assumed $P$ was compact to
apply Theorem \ref{slevthm}. This assumption is in fact unnecessary. 
As in the discussion after Theorem \ref{slevthm}, for any $p\in P$
the $\phi_t$ exist near $p$ for $t\in(-\ep,\ep)$ and some $\ep>0$. We 
can then use $G$-equivariance to extend $\phi_t$ uniquely to all of~$P$.
\end{proof}

This means that special Lagrangian submanifolds of cohomogeneity 
one in $\C^m$ are relatively easy to construct and classify. The
strategy is to first identify all the suitable Lie subgroups $G$ 
in $\SU(m)\lt\C^m$ which admit moment maps. Note that though not 
all subgroups of $\SU(m)\lt\C^m$ admit moment maps, the symmetry
group $\Sym(N)$ of a Lagrangian submanifold always admits a 
moment map, and so subgroups without moment maps are excluded 
as they cannot preserve any special Lagrangian submanifolds.

Once we have chosen a Lie subgroup $G$ with moment map $\mu$, we 
then work out the types of $G$-orbit $\cal O$ in $\mu^{-1}(c)$ for 
$c\in Z(\g^*)$, and see if any have dimension $m-1$. Clearly we 
must have $\dim G\ge m-1$ to get any suitable orbits. But if 
$\dim G$ is too large then there won't be any suitable orbits either.

We then find the cohomogeneity one $G$-invariant special 
Lagrangian $m$-folds in $\C^m$ by solving a first-order ordinary 
differential equation in $(m-1)$-dimensional $G$-orbits. This can 
often be done explicitly, or failing this, a qualitative description 
of the solutions can be given.

Theorem \ref{sy4thm} says only that this o.d.e.\ is soluble
for small $t$. Solutions generally exist in some open interval
$I$ in $\R$. As $t$ approaches the ends of the interval, two 
things can happen: the orbit $\phi_t({\cal O})$ can go off to 
infinity, or it can collapse down to another $G$-orbit of smaller
dimension. By including this $G$-orbit in $N$, one sometimes gets
a closed, nonsingular submanifold in~$\C^m$.

\section{Examples of cohomogeneity one SL $m$-folds}
\label{sy5}

We now give some examples of cohomogeneity one special Lagrangian
$m$-folds in $\C^m$, as in \S\ref{sy41}. We begin with three examples in 
$\C^3$, the first taken from Harvey and Lawson~\cite[\S III.3.A]{HaLa}.

\begin{ex} Let $G\cong T^2$ be the group of diagonal matrices 
in $\SU(3)$, so that each $\ga\in G$ acts on $\C^3$ by 
\begin{equation*}
\ga:(z_1,z_2,z_3)\mapsto({\rm e}^{i\th_1}z_1,{\rm e}^{i\th_2}z_2,
{\rm e}^{i\th_3}z_3)
\end{equation*}
for some $\th_1,\th_2,\th_3\in\R$ with $\th_1+\th_2+\th_3=0$. 
As $G$ is abelian, $Z(\g^*)=\g^*$, and for any $c\in\g^*$ the
generic $G$-orbit in $\mu^{-1}(c)$ is a copy of $T^2$, and
so 2-dimensional.

Thus we can apply Theorem \ref{sy4thm} to find a family of 
$T^2$-invariant special Lagrangian 3-folds in $\C^3$, by
solving an ordinary differential equation. Let $a_1,a_2$ and $b$ 
be real numbers, and define a subset $L_{a_1,a_2,b}$ in $\C^3$ by
\begin{align*}
L_{a_1,a_2,b}=\bigl\{(z_1,&z_2,z_3)\in\C^3:
\ms{z_1}-\ms{z_3}=a_1,\\
&\ms{z_2}-\ms{z_3}=a_2,\quad
\Im(z_1z_2z_3)=b\bigr\}.
\end{align*}
Harvey and Lawson show that $L_{a_1,a_2,b}$ is a $T^2$-invariant 
SL 3-fold in $\C^3$, and it's also easy to see 
that any connected, $T^2$-invariant SL 3-fold is 
a subset of some~$L_{a_1,a_2,b}$.

Here the equations $\ms{z_1}-\ms{z_3}=a_1$ and $\ms{z_2}-\ms{z_3}=a_2$
are the moment maps of $G$, which must be constant on $L_{a_1,a_2,b}$
by Proposition \ref{sy4prop2}. The third equation $\Im(z_1z_2z_3)=b$
can also be interpreted as a kind of generalized moment map equation,
associated to the 3-form~$\Im\Om$.

The $L_{a_1,a_2,b}$ are not all nonsingular. In fact one can show:
\begin{itemize}
\setlength{\parsep}{0pt}
\setlength{\itemsep}{0pt}
\item[(i)] $L_{0,0,0}$ has an isolated singular point at $(0,0,0)$.
\item[(ii)] $L_{r^2,0,0}$ is singular on the circle 
$\bigl\{(z_1,0,0):\md{z_1}=r\bigr\}$ for~$r>0$.
\item[(iii)] $L_{0,r^2,0}$ is singular on the circle
$\bigl\{(0,z_2,0):\md{z_2}=r\bigr\}$ for~$r>0$.
\item[(iv)] $L_{-r^2,-r^2,0}$ is singular on the circle
$\bigl\{(0,0,z_3):\md{z_3}=r\bigr\}$ for~$r>0$.
\item[(v)] All other $L_{a_1,a_2,b}$ are nonsingular, 
and diffeomorphic to~$\R\t T^2$.
\end{itemize}
\label{c3ex1}
\end{ex}

Observe that $\C^3$ is {\it fibred} by this family of special 
Lagrangian 3-folds; that is, there is exactly one passing through
each point in $\C^3$. Here is a second example, adapted from Harvey 
and Lawson \cite[\S III.3.B]{HaLa}. It can also be obtained by
applying Theorem \ref{sy6thm2} to the special Lagrangian cone
$\R^3$ in~$\C^3$.

\begin{ex} Let $G$ be the subgroup $\SO(3)$ of $\SU(3)$,
embedded as $3\t 3$ real matrices in the obvious way.
Then the moment map of $G$ is
\begin{equation*}
\mu(z_1,z_2,z_3)=\bigl(\Im(z_1\bar z_2),\Im(z_2\bar z_3),
\Im(z_3\bar z_1)\bigr).
\end{equation*}
As $Z(\g^*)=\{0\}$, any $G$-invariant special Lagrangian 3-fold 
lies in~$\mu^{-1}(0)$. 

Now all points in $\mu^{-1}(0)$ may be written as 
$(\la x_1,\la x_2,\la x_3)$, where $\la\in\C$ and 
$x_1,x_2,x_3$ are {\it real}, and normalized so that 
$x_1^2+x_2^2+x_3^2=1$. Therefore the orbits of $G$ in 
$\mu^{-1}(0)$ are ${\cal O}_\la$ for $\la\in\C$, where
\begin{equation*}
{\cal O}_\la=\bigl\{(\la x_1,\la x_2,\la x_3):
x_j\in\R,\quad x_1^2+x_2^2+x_3^2=1\bigr\}.
\end{equation*}
Clearly ${\cal O}_0$ is a point, and ${\cal O}_\la\cong{\cal S}^2$
if $\la\ne 0$, and~${\cal O}_\la={\cal O}_{-\la}$.

We can therefore interpret the o.d.e.\ on $G$-orbits in $\mu^{-1}(0)$
discussed in Theorem \ref{sy4thm} as an o.d.e.\ on $\la$. Calculation
shows that it is $\d\la/\d t=\bar\la^2$, where $\la=\la(t)$. 
Hence we see that $\d(\la^3)/\d t=3\md{\la}^4$, which is real, and 
so $\d(\Im(\la^3))/\d t=0$. Thus the integral curves of the o.d.e.\ 
are of the form $\Im(\la^3)=c$ for~$c\in\R$.

So for each $c\in\R$, define
\begin{equation*}
N_c=\bigl\{(\la x_1,\la x_2,\la x_3):
x_j\in\R,\quad x_1^2+x_2^2+x_3^2=1,\quad
\la\in\C,\quad \Im(\la^3)=c\bigr\}.
\end{equation*}
Then $N_c$ is a special Lagrangian 3-fold in $\C^3$. When $c=0$,
it is a singular union of three copies of $\R^3$ intersecting
at 0, and when $c\ne 0$ it is nonsingular, the disjoint union
of three copies of $\R\t{\cal S}^2$. Note also that~$N_c=N_{-c}$.
\label{c3ex2}
\end{ex}

Here is a rather trivial example.

\begin{ex} Let $G$ be the subgroup $U(1)\t\R$ of $\SU(3)\ltimes\C^3$,
acting by
\begin{equation*}
(e^{i\th},x):(z_1,z_2,z_3)\longmapsto(e^{i\th}z_1,e^{-i\th}z_2,x+z_3),
\end{equation*}
for $\th\in[0,2\pi)$ and $x\in\R$. The moment map of $G$ is
\begin{equation*}
\mu(z_1,z_2,z_3)=\bigl(\ms{z_1}-\ms{z_2},\Im z_3\bigr),
\end{equation*}
and $\Z(\g^*)=\g^*=\R^2$, and $G$-orbits are copies of
${\cal S}^1\t\R$ unless $z_1=z_2=0$, when they are copies of~$\R$.

As in Example \ref{c3ex1}, we find the following. Let $a,b,c\in\R$, 
and define
\begin{equation*}
N_{a,b,c}=\bigl\{(z_1,z_2,z_3)\in\C^3:\ms{z_1}-\ms{z_2}=a,\quad
\Re(z_1z_2)=b,\quad \Im(z_3)=c\bigr\}.
\end{equation*}
Then $N_{a,b,c}$ is a $G$-invariant special Lagrangian 3-fold in 
$\C^3$. If $a=b=0$ then $N_{a,b,c}$ is the singular union of two
copies of $\R^3$ intersecting in $\R$, and otherwise $N_{a,b,c}$
is nonsingular and diffeomorphic to~${\cal S}^1\t\R^2$. 
\label{c3ex3}
\end{ex}

Again, $\C^3$ is fibred by these $N_{a,b,c}$. Note that 
$N_{a,b,c}$ is actually the product of lower-dimensional 
SL submanifolds in $\C^2$ and $\C$. It turns out that these 
three examples represent all cohomogeneity one SL
3-folds in~$\C^3$.

\begin{thm} Every homogeneous special Lagrangian $3$-fold in 
$\C^3$ is conjugate under $\SU(3)\ltimes\C^3$ to the standard\/
$3$-plane $\R^3\subset\C^3$. Every cohomogeneity one special 
Lagrangian $3$-fold in $\C^3$ is conjugate under $\SU(3)\ltimes\C^3$ 
to a subset of\/ $\R^3$ or one of the $3$-folds of Examples 
\ref{c3ex1}, \ref{c3ex2} and\/~\ref{c3ex3}.
\label{sy5thm}
\end{thm}

We leave the proof as an exercise; one uses the classification of 
Lie groups to identify all Lie subgroups $G$ of $\SU(3)\ltimes\C^3$,
and show that either $G$ has no suitable 2-dimensional orbits in
$\mu^{-1}(c)$ for $c\in Z(\g^*)$, or else that the $G$-invariant
submanifolds reduce to one of the cases in the theorem. Note that
$\Sym(\R^3)$ is $\SO(3)\ltimes\R^3$, which has several Lie subgroups 
$G$ leading to subsets of~$\R^3$.

Next we give some higher-dimensional examples. Here is an example 
generalizing Example \ref{c3ex1}, taken from~\cite[\S III.3.A]{HaLa}.

\begin{ex} Let $G\cong T^{m-1}$ be the group of diagonal matrices 
in $\SU(m)$, so that each $\ga\in G$ acts on $\C^m$ by 
\begin{equation*}
\ga:(z_1,\ldots,z_m)\mapsto({\rm e}^{i\th_1}z_1,\ldots,
{\rm e}^{i\th_m}z_m)
\end{equation*}
for some $\th_1,\ldots,\th_m\in\R$ with $\th_1+\cdots+\th_m=0$. 

Let $a_1,\ldots,a_{m-1}$ and $b$ be real numbers. Define 
$L_{a_1,\ldots,a_{m-1},b}$ by
\begin{align*}
L_{a_1,\ldots,a_{m-1},b}=\Bigl\{(z_1,&\ldots,z_m)\in\C^m:
\text{$\ms{z_j}-\ms{z_m}=a_j$ for $j=1,\ldots,m-1$,}\\
&\text{and}\quad
\begin{cases}\Re(z_1\ldots z_m)=b & \text{if $m$ is even}\\
\Im(z_1\ldots z_m)=b & \text{if $m$ is odd}\end{cases}\,\,
\Bigr\}.
\end{align*}
Then Harvey and Lawson show that $L_{a_1,\ldots,a_{m-1},b}$
is a special Lagrangian $m$-fold in $\C^m$. If $b\ne 0$ then
$L_{a_1,\ldots,a_{m-1},b}$ is nonsingular and diffeomorphic
to $T^{m-1}\t\R$. When $b=0$, it may be nonsingular or have
various kinds of singularity, depending on the values 
of~$a_1,\ldots,a_{m-1}$.
\label{cmex1}
\end{ex}

Again, $\C^m$ is fibred by these special Lagrangian 
$m$-folds. Here is another example of Harvey and Lawson 
\cite[\S III.3.B]{HaLa}, generalizing Example \ref{c3ex2}.
It can also be derived from Theorem \ref{sy6thm2} with
$C=\R^m$ in~$\C^m$.

\begin{ex} Let $G$ be $\SO(m)$ in $\SU(m)$. For each $c\in\R$, define
\begin{equation*}
N_c=\bigl\{(\la x_1,\ldots,\la x_m):
x_j\in\R,\quad x_1^2+\cdots+x_m^2=1,\quad
\la\in\C,\quad \Im(\la^m)=c\bigr\}.
\end{equation*}
Then $N_c$ is a special Lagrangian $m$-fold in $\C^m$. When $c=0$,
it is a singular union of $m$ copies of $\R^m$ intersecting at 0.
If $m$ is even and $c\ne 0$, $N_c$ is a nonsingular, disjoint
union of $m/2$ copies of $\R\t{\cal S}^{m-1}$. If $m$ is odd 
and $c\ne 0$, $N_c$ is a nonsingular, disjoint union of $m$ 
copies of $\R\t{\cal S}^{m-1}$, and~$N_c=N_{-c}$.
\label{cmex2}
\end{ex}

We conclude with a more complex example due to 
Marshall~\cite[\S 3.4]{Mars}.

\begin{ex} The usual action of $\SU(2)$ on $\C^2$ induces an 
action of $\SU(2)$ on $S^3\C^2$. Identifying $S^3\C^2$ with 
$\C^4$ in an appropriate way, this defines a subgroup $G$
of $\SU(4)$ isomorphic to $\SU(2)$. Calculation shows that
we may take the Lie algebra $\g$ of $G$ to be spanned by
\begin{equation*}
\begin{pmatrix}
3i & 0 & 0 & 0 \\ 0 & i & 0 & 0 \\ 0 & 0 &-i & 0 \\ 0 & 0 & 0 &-3i
\end{pmatrix},\;
\begin{pmatrix}
0 & \sqrt{3} & 0 & 0 \\ -\sqrt{3} & 0 & 2 & 0 \\
0 & -2 & 0 & \sqrt{3} \\ 0 & 0 & -\sqrt{3} & 0 
\end{pmatrix},\;
\begin{pmatrix}
0 & i\sqrt{3} & 0 & 0 \\ i\sqrt{3} & 0 & 2i & 0 \\
0 & 2i & 0 & i\sqrt{3} \\ 0 & 0 & i\sqrt{3} & 0 
\end{pmatrix}.
\end{equation*}
For each $d\in\R$, define $N_d\subset\C^4$ by
\begin{align*}
N_d=\Bigl\{(z_1,\ldots,z_4)\in\C^4:
&\,\sqrt{3}(z_1\bar z_2+z_3\bar z_4)+2z_2\bar z_3=0,\\
&\,3\ms{z_1}+\ms{z_2}-\ms{z_3}-3\ms{z_4}=0,\\
\Im\Bigl(2\sqrt{3}(z_1z_3^2+&z_2^3z_4)
-9z_1z_2z_3z_4+{9\over 2}z_1^2z_4^2
-{3\over 2}z_2^2z_3^2\Bigr)=d\Bigr\}.
\end{align*}
Then $N_d$ is a $G$-invariant SL 4-fold in $\C^4$. Here the first two 
equations defining $N_d$ say that the moment map $\mu$ is zero. 

It can be shown that $N_0$ is the union of two cones on ${\cal S}^3/\Z_3$,
with one singular point at 0, where if we identify ${\cal S}^3$ with the 
unit sphere in $\C^2$, then the $\Z_3$-action is generated by $(w_1,w_2)
\mapsto({\rm e}^{2\pi i/3}w_1,{\rm e}^{2\pi i/3}w_2\bigr)$. Similarly,
if $d\ne 0$ then $N_d$ is nonsingular and diffeomorphic 
to~$\R\t{\cal S}^3/\Z_3$.
\label{c4ex1}
\end{ex}             

\section{Special Lagrangian cones in $\C^m$}
\label{sy6}

In \S\ref{sy4} we studied symmetries of special Lagrangian
$m$-folds in $\C^m$, which were required to preserve the metric
$g$, K\"ahler form $\om$ and complex volume form $\Om$ on $\C^m$.
However, there is a more general kind of automorphism of $\C^m$ 
which does not preserve $g$, $\om$ and $\Om$, but does 
preserve the idea of SL submanifolds in~$\C^m$.

\begin{dfn} Let $\R_+$ be the group of positive real numbers 
under multiplication, and let $t\in\R_+$ act on $\C^m$ by
$(z_1,\ldots,z_m)\,{\buildrel t\over\longmapsto}\,(tz_1,\ldots,tz_m)$.
We call this action of $t$ on $\C^m$ a {\it dilation}. Let $N$ be a 
real $m$-dimensional submanifold of $\C^m$. Then $tN$ is also a real 
$m$-dimensional submanifold of $\C^m$ for each $t\in\R_+$. Clearly,
$N$ is special Lagrangian if and only if $tN$ is, so that dilations 
preserve the idea of SL submanifolds in~$\C^m$. 

Combining dilations with the automorphisms $\SU(m)\lt\C^m$ of $\C^m$ 
gives a group $\bigl(\R_+\!\t\SU(m)\bigr)\lt\C^m$ acting on $\C^m$ 
preserving SL submanifolds. If $N$ is a special 
Lagrangian submanifold of $\C^m$, define the {\it generalized 
symmetry group} $\GSym(N)$ of $N$ to be the Lie subgroup $G$ of 
$\bigl(\R_+\!\t\SU(m)\bigr)\lt\C^m$ preserving $N$, and define
the {\it restricted generalized symmetry group} $\GSym^0(N)$ to 
be the identity component of~$\GSym(N)$.
\label{gsymndef}
\end{dfn}

The symmetry group $\Sym(N)$ of \S\ref{sy4} is a normal subgroup 
of $\GSym(N)$, and $\GSym(N)/\Sym(N)$ is a subgroup of $\R_+$. 
It is convenient to distinguish between $\Sym(N)$ and $\GSym(N)$,
because $\Sym^0(N)$ always admits a moment map, but in general
$\GSym^0(N)$ has no moment map as it doesn't preserve $\om$.

Submanifolds invariant under dilations are called {\it cones}.

\begin{dfn} Let $C$ be a (singular) submanifold of $\C^m$. We call 
$N$ a {\it cone} in $\C^m$, with vertex 0, if $0\in C$ and $tC=C$ 
for all $t\in\R_+$. That is, $C$ is a cone if it is invariant under 
dilations. Let $C$ be a cone, and define $\Si=\bigl\{{\bf z}\in C:
\md{{\bf z}}=1\bigr\}=C\cap{\cal S}^{2m-1}$. Then $C=\bigl\{t{\bf z}:
{\bf z}\in\Si$, $t\in[0,\iy)\bigr\}$. We call $\Si$ the {\it link} 
of $C$, and $C$ the {\it cone on}~$\Si$.

A closed submanifold $N$ in $\C^m$ is called {\it Asymptotically 
Conical}, or {\it AC} for short, if there exists a closed cone 
$N_0$ in $\C^m$ with isolated singular point at 0, such that $N$ 
is asymptotic to $N_0$ to order $O(r^{-1})$ as $r\ra\iy$, where 
$r$ is the radius function in $\C^m$. We call $N_0$ the 
{\it asymptotic cone} of~$N$.
\label{conedef}
\end{dfn}

If $N$ is a special Lagrangian cone, then $\GSym(N)=\R_+\t\Sym(N)$.
We will be particularly interested in conical and asymptotically
conical SL submanifolds in $\C^m$, as they can 
provide local models for singularities of SL 
$m$-folds in Calabi--Yau $m$-folds. Here is why. Let $N$ be a 
nonsingular AC special Lagrangian $m$-fold in $\C^m$, asymptotic
to a singular cone $N_0$. If $t>0$ then $tN$ is also AC, and 
$tN\ra N_0$ as $t\ra 0_+$. 

Thus the singular SL submanifold $N_0$ is the 
limit of the family of nonsingular SL submanifolds 
$\{tN:t>0\}$. So AC special Lagrangian submanifolds provide local 
models for how singularities can develop in families of nonsingular 
SL submanifolds.

\subsection{Special Lagrangian cones of cohomogeneity one}
\label{sy61}

In \S\ref{sy41} we studied special Lagrangian $m$-folds $N$ in 
$\C^m$ upon which $\Sym(N)$ acts with cohomogeneity one. We can 
also consider $N$ upon which the {\it generalized}\/ symmetry 
group $\GSym(N)$ above acts with cohomogeneity one. In particular, 
if $N$ is a cone then $\GSym(N)=\R_+\!\t\Sym(N)$, so that $\GSym(N)$ 
has cohomogeneity one when $\Sym(N)$ has cohomogeneity two. Here is 
a result on cohomogeneity one cones, similar to Theorem~\ref{sy4thm}.

\begin{thm} Let\/ $G$ be a connected Lie subgroup of\/ 
$\SU(m)$ with Lie algebra $\g$ and moment map $\mu:\C^m\ra\g^*$, 
with\/ $\mu(0)=0$. Let\/ $\cal O$ be an oriented orbit of\/ $\R_+\!\t G$ 
in $\C^m$ with\/ $\dim{\cal O}=m-1$, and suppose ${\cal O}\subset\mu^{-1}(0)$. 
Then there exists a locally unique, $\R_+\!\t G$-invariant special 
Lagrangian cone $N$ in $\C^m$ containing $\cal O$. Furthermore 
$N\subset\mu^{-1}(0)$, and\/ $N$ is fibred by $\R_+\!\t G$-orbits 
isomorphic to $\cal O$ near $\cal O$. Thus $N$ is locally diffeomorphic 
to $(-\ep,\ep)\t{\cal O}$, for some $\ep>0$, and we can think of\/ $N$ 
as a smooth curve of\/ $\R_+\!\t G$-orbits.
\label{sy6thm1}
\end{thm}

The proof is very similar to Theorem \ref{sy4thm}, so we will
not give it, but here are a few comments. Since we are constructing
cones in $\C^m$ we exclude translations from our symmetry group, 
as they would move the vertex. Thus we take $G\subset\SU(m)$ rather 
than $G\subset\SU(m)\lt\C^m$. The moment map $\mu$ for $G$ always 
exists, and we specify it uniquely by requiring that~$\mu(0)=0$. 

Suppose $N$ is an $(\R_+\!\t G)$-invariant special Lagrangian cone. 
Then $N\subset\mu^{-1}(c)$ for $c\in Z(\g^*)$, by Proposition 
\ref{sy4prop2}. But $\mu(t{\bf z})=t^2\mu({\bf z})$ for $t\in\R_+$,
which forces $c=t^2c$ as $N=tN$, so $c=0$. Thus $N\subset\mu^{-1}(0)$.
Conversely, following Proposition \ref{sy4prop1} one can show that if 
$\cal O$ is an $\R_+\!\t G$-orbit then $\om\vert_{\cal O}\equiv 0$ if
and only if ${\cal O}\subset\mu^{-1}(0)$, and this is used to apply
Theorems \ref{hlthm} and~\ref{slevthm}.

\subsection{AC special Lagrangian $m$-folds from cones}
\label{sy62}

We shall now show that given {\it any} special Lagrangian cone 
in $\C^m$, one can automatically construct a 1-parameter family
of asymptotically conical SL $m$-folds in $\C^m$ from it. This was 
first noticed by Castro and Urbano \cite[Remark 1, p.~81--2]{CaUr2}
in the context of minimal Lagrangian submanifolds, proved independently
by Haskins \cite[Th.~A]{Hask}, and also independently by the author.

\begin{thm} Let\/ $C$ be a closed special Lagrangian cone in $\C^m$
with isolated singular point\/ $0$, set\/ $\Si=\bigl\{{\bf z}\in C:
\md{{\bf z}}=1\bigr\}$, and for each $c>0$ define
\e
\begin{split}
N_c&=\bigl\{\la{\bf z}:{\bf z}\in\Si,\;\> \la\in\C,\;\>
\Im(\la^m)=c^m,\;\> \arg(\la)\in(0,\pi/m)\bigr\}\\
&=\bigl\{c(\sin(m\th))^{-1/m}{\rm e}^{i\th}{\bf z}:{\bf z}\in\Si,\quad
\th\in(0,\pi/m)\bigr\}.
\end{split}
\label{sy6eq1}
\e
Then $N_c$ is an immersed AC special Lagrangian $m$-fold 
in $\C^m$ diffeomorphic to $\Si\t\R$, and asymptotic 
to~$C\cup {\rm e}^{i\pi/m}C$.
\label{sy6thm2}
\end{thm}

\begin{proof} One can prove this quite simply by relating
the tangent spaces $T_pN_c$ of $N_c$ to those of $\Si$, and 
showing that each $T_pN_c$ is special Lagrangian. But we will 
instead give a proof using Theorem~\ref{slevthm}.

Now $\Si$ is a compact, nonsingular Riemannian $(m-1)$-manifold,
with a natural orientation. Let $\chi$ be the unique positive
section of $\La^{m-1}T\Si$ with $\md{\chi}\equiv 1$. It follows 
easily from Theorem \ref{morrthm} that $\Si$ and $\chi$ are
real analytic. Consider a 1-parameter family $\bigl\{\phi_t:
t\in(-\ep,\ep)\bigr\}$ of maps $\phi_t:\Si\ra\C^m$ defined by 
$\phi_t:{\bf z}\mapsto\la(t){\bf z}$, where $\la(t):(-\ep,\ep)
\ra\C\sm\{0\}$ is a differentiable function.

Calculation shows that the $\phi_t$ satisfy the evolution equation
\eq{dphiteq} of Theorem \ref{slevthm} if and only if $\la$ satisfies
the o.d.e.
\e
{\d\la\over\d t}=\bar\la^{m-1}.
\label{sy6eq2}
\e
But then $\d(\la^m)/\d t=m\la^{m-1}\bar\la^{m-1}=m\md{\la}^{2m-2}$,
which is real. So $\Im(\la^m)$ is constant along the integral curves
of \eq{sy6eq2}. In particular, if $c>0$ then $\bigl\{\la\in\C:
\Im(\la^m)=c^m$, $\arg(\la)\in(0,\pi/m)\bigr\}$ is an integral 
curve of \eq{sy6eq2}, and the result then follows quickly from 
Theorem~\ref{slevthm}.
\end{proof}

\section{$\U(1)^{m-2}$-invariant cones in $\C^m$}
\label{sy7}

We shall now apply Theorem \ref{sy6thm1} to the example of 
$G=\U(1)^{m-2}$ in $\SU(m)$. The case $m=3$ has already been 
analyzed by Mark Haskins \cite[\S 3--\S 5]{Hask}, using 
somewhat different techniques.

Let $m\ge 3$ and $a_1\le\cdots\le a_m$ be integers, not all zero, with 
highest common factor 1, and with $a_1+\cdots+a_m=0$. Note that this
implies that $a_1<0$ and $a_m>0$. Define a subgroup $G\subset\U(1)^m$ 
to be
\begin{equation*}
\bigl\{({\rm e}^{i\al_1},\ldots,{\rm e}^{i\al_m})\in\U(1)^m:
\al_j\in\R,\;\> \al_1\!+\!\cdots\!+\!\al_m\!=\!0,\;\> 
a_1\al_1\!+\!\cdots\!+\!a_m\al_m\!=\!0\bigr\}.
\end{equation*}
Then $G\cong\U(1)^{m-2}$, because the $a_j$ are integers. If we
instead allowed the $a_j$ to be real numbers then $G$ would be
isomorphic to $\U(1)^k\t\R^{m-2-k}$ for some $0\le k\le m-2$.
Most of the analysis below would still work, but the resulting
SL $m$-folds would not be as interesting.

Let $G$ act on $\C^m$ in the obvious way, by
\begin{equation*}
({\rm e}^{i\al_1},\ldots,{\rm e}^{i\al_m}):(z_1,\ldots,z_m)
\mapsto({\rm e}^{i\al_1}z_1,\ldots,{\rm e}^{i\al_m}z_m).
\end{equation*}
Then $G\subset\SU(m)$. We shall apply Theorem \ref{sy6thm1} to construct 
$\R_+\!\t G$-invariant SL $m$-folds $N$ in $\C^m$, which 
we will regard as the total space of a 1-parameter family of 
$\R_+\!\t G$-orbits ${\cal O}_t$ for~$t\in(-\ep,\ep)$. 

So, for $t\in(-\ep,\ep)$ define $\phi_t:\R_+\!\t G\ra\C^m$ by
\e
\phi_t:\bigl(r,({\rm e}^{i\al_1},\ldots,{\rm e}^{i\al_m})\bigr)\mapsto
\bigl(r{\rm e}^{i\al_1}w_1(t),\ldots,r{\rm e}^{i\al_m}w_m(t)\bigr)
\label{sy7eq1}
\e
for $r\in \R_+$ and $({\rm e}^{i\al_1},\ldots,{\rm e}^{i\al_m})\in G$,
where $w_1,\ldots,w_m:(-\ep,\ep)\ra\C$ are differentiable functions.
Let ${\cal O}_t=\Im\phi_t$, and let $N$ be the union of the ${\cal O}_t$
for $t\in(-\ep,\ep)$. We will shortly use equation \eq{dphiteq} of Theorem 
\ref{slevthm} to derive an {\it evolution equation} for $\phi_t$
and $w_1,\ldots,w_m$, which will imply that $N$ is special Lagrangian
in $\C^m$, with phase~$i^{m-2}$.

First, we consider what the conditions that ${\cal O}\subset\mu^{-1}(0)$ 
and $\dim{\cal O}=m-1$ in Theorem \ref{sy6thm1} mean for $w_1,\ldots,w_m$.
The Lie algebra of $\U(1)^m$ is $\R^m$, and the Lie algebra $\g$ of 
$G$ is the subspace of $(x_1,\ldots,x_m)\in\R^m$ such that 
$x_1+\cdots+x_m=0$ and $a_1x_1+\cdots+a_mx_m=0$. Thus, the condition 
that $(w_1,\ldots,w_m)\in\mu^{-1}(0)$ is that $x_1\ms{w_1}+\cdots
+x_m\ms{w_m}=0$ whenever $x_1,\ldots,x_m\in\R$, $x_1+\cdots+x_m=0$ 
and~$a_1x_1+\cdots+a_mx_m=0$. 

Clearly, this is true if and only if
\e
\ms{w_j}=a_ju+v,\quad j=1,\ldots,m,
\quad\text{for some $u,v\in\R$.}
\label{sy7eq2}
\e
This is also the condition that ${\cal O}_t\subset\mu^{-1}(0)$,
and by the discussion after Theorem \ref{sy6thm1} this is
equivalent to $\om\vert_{{\cal O}_t}\equiv 0$. Also, $\dim{\cal O}_t=m-1$ 
holds if no more than one of $w_1(t),\ldots,w_m(t)$ are zero 
for all $t\in(-\ep,\ep)$. For simplicity we make the stronger
assumption that $w_j(t)$ is nonzero for all $j$ and~$t$.

Next, we use equation \eq{dphiteq} to derive an o.d.e.\ for $\phi_t$
and $w_1,\ldots,w_m$. To do this we need an $\R_+\!\t G$-invariant
section $\chi$ of $\La^{m-1}T(\R_+\!\t G)$. In terms of the natural
local coordinates $r$ and $\al_1,\ldots,\al_m$ on $\R_+\t\U(1)^m$,
calculation shows that we may take $\chi$ to be
\begin{equation*}
\chi={2r\over m}{\pd\over\pd r}\w\!\!\!\!\!\!\sum_{1\le j<k\le m}\!\!\!\!\!
(-1)^{j\!+\!k\!-\!1}(a_j-a_k)\pd_1\!\w\!\cdots
\pd_{j\!-\!1}\!\w\!\pd_{j\!+\!1}\!\w\!\cdots\!\w\!
\pd_{k\!-\!1}\!\w\!\pd_{k\!+\!1}\!\w\!\cdots\!\w\!\pd_m,
\end{equation*}
where $\pd_j=\pd/\pd\al_j$. To apply \eq{dphiteq}, we need an 
expression for $(\phi_t)_*(\chi)$ at a point $(z_1,\ldots,z_m)$ 
in ${\cal O}_t$. Using the equations
\begin{align*}
(\phi_t)_*\Bigl(r{\pd\over\pd r}\Bigr)&=
z_1{\pd\over\pd z_1}+\cdots+z_m{\pd\over\pd z_m}+
\bar z_1{\pd\over\pd\bar z_1}+\cdots+\bar z_m{\pd\over\pd\bar z_m}\\
\text{and}\quad (\phi_t)_*(\pd_j)&=
iz_j{\pd\over\pd z_j}-i\bar z_j{\pd\over\pd\bar z_j}
\end{align*}
we can do this, and the result is rather complicated. 

But all we will actually need is the $(m-1,0)$ component of 
$(\phi_t)_*(\chi)$, throwing away all terms in $\pd/\pd\bar z_j$,
and calculation shows that this is given by
\begin{align*}
(\phi_t)_*(\chi)^{(m-1,0)}=
2i^{m-2}\sum_{j=1}^m
&(-1)^{j-1}a_jz_1\cdots z_{j\!-\!1}z_{j\!+\!1}\cdots z_m \cdot \\
&{\pd\over\pd z_1}\!\w\cdots\w\!{\pd\over\pd z_{j\!-\!1}}\!\w\! 
{\pd\over\pd z_{j\!+\!1}}\!\w\cdots\w\!{\pd\over\pd z_m}.
\end{align*}
As $\Om$ is an $(m,0)$-tensor, we see that the contraction
of $(\phi_t)_*(\chi)$ with $\Om$ is the same as that of 
$(\phi_t)_*(\chi)^{(m-1,0)}$ with $\Om$. Hence, using the index 
notation for tensors on $\C^m$, we get
\begin{equation*}
(\phi_t)_*(\chi)^{b_1\ldots b_{m-1}}\Om_{b_1\ldots b_{m-1}b_m}=
2i^{m-2}\sum_{j=1}^m
a_jz_1\cdots z_{j\!-\!1}z_{j\!+\!1}\cdots z_m(\d z_j)_{b_m}.
\end{equation*}
Multiplying by $(-i)^{m-2}$ and contracting with 
$g^{b_mc}$ gives
\begin{gather*}
(\phi_t)_*(\chi)^{b_1\ldots b_{m-1}}
\bigl((-i)^{m-2}\Om\bigr)_{b_1\ldots b_{m-1}b_m}g^{b_mc}=\\
2\sum_{j=1}^ma_jz_1\cdots z_{j\!-\!1}z_{j\!+\!1}\cdots z_m 
\Bigl({\pd\over\pd\bar z_j}\Bigr)^c.
\end{gather*}

Since $(\phi_t)_*(\chi)$ and $g$ are real tensors, taking
real parts gives
\begin{gather*}
(\phi_t)_*(\chi)^{b_1\ldots b_{m-1}}
\Re\bigl((-i)^{m-2}\Om\bigr)_{b_1\ldots b_{m-1}b_m}g^{b_mc}=\\
\sum_{j=1}^ma_j\overline{z_1\cdots z_{j\!-\!1}z_{j\!+\!1}\cdots z_m} 
\Bigl({\pd\over\pd z_j}\Bigr)^c+
\sum_{j=1}^ma_jz_1\cdots z_{j\!-\!1}z_{j\!+\!1}\cdots z_m 
\Bigl({\pd\over\pd\bar z_j}\Bigr)^c.
\end{gather*}
Now by Theorem \ref{slevthm}, a sufficient condition for $N$ 
to be special Lagrangian with phase $i^{m-2}$ is 
\begin{equation*}
\left({\d\phi_t\over\d t}\right)^c=(\phi_t)_*(\chi)^{b_1\ldots b_{m-1}}
\bigl(\Re(-i)^{m-2}\Om\bigr)_{b_1\ldots b_{m-1}b_m}g^{b_mc},
\end{equation*}
where we have inserted the factor $(-i)^{m-2}$ in front of $\Om$
to get $N$ of phase $i^{m-2}$ rather than~1.

Setting $(z_1,\ldots,z_m)$ to be $\phi_t(1,1,\ldots,1)=(w_1,\ldots,w_m)$
and combining the last two equations, we find that a sufficient
condition for $N$ to be special Lagrangian with phase $i^{m-2}$ is 
\begin{gather*}
\sum_{j=1}^m{\d w_j\over\d t}\Bigl({\pd\over\pd z_j}\Bigr)^c
+\sum_{j=1}^m{\d\bar w_j\over\d t}\Bigl({\pd\over\pd\bar z_j}\Bigr)^c=\\
\sum_{j=1}^ma_j\,\overline{w_1\cdots w_{j\!-\!1}w_{j\!+\!1}\cdots w_m} 
\Bigl({\pd\over\pd z_j}\Bigr)^c+
\sum_{j=1}^ma_j\,w_1\cdots w_{j\!-\!1}w_{j\!+\!1}\cdots w_m 
\Bigl({\pd\over\pd\bar z_j}\Bigr)^c.
\end{gather*}
Equating coefficients, this is true if
\e
{\d w_j\over\d t}=a_j\,\overline{w_1\cdots w_{j\!-\!1}w_{j\!+\!1}\cdots w_m}
\quad\text{for $j=1,\ldots,m$.}
\label{sy7eq3}
\e

Now we assumed above that \eq{sy7eq2} holds for $w_1,\ldots,w_m$, so
we should check that our evolution equation \eq{sy7eq3} preserves
$w_j$ of this form. From \eq{sy7eq3} we get
\begin{equation*}
{\d\ms{w_j}\over\d t}=w_j{\d\bar w_j\over\d t}+\bar w_j{\d w_j\over\d t}=
a_jw_1\cdots w_m+a_j\ov{w_1\cdots w_m}=2a_j\Re(w_1\cdots w_m).
\end{equation*}
Thus the evolution \eq{sy7eq3} does preserve $w_1,\ldots,w_m$ of
the form \eq{sy7eq2}, and $u,v$ are functions of $t$ satisfying
\begin{equation*}
{\d u\over\d t}=2\Re(w_1\cdots w_m), \qquad {\d v\over\d t}=0. 
\end{equation*}
So $v$ is constant. Since $a_1+\cdots+a_m=0$ we have
$\ms{w_1}+\cdots+\ms{w_m}=mv$, and $v>0$. As rescaling
all the $z_j$ by a positive constant leads to the same
SL $m$-fold, we may as well fix~$v=1$.

We summarize our progress so far in the following theorem.

\begin{thm} Suppose $w_1,\ldots,w_m:(-\ep,\ep)\ra\C\sm\{0\}$ and\/ 
$u:(-\ep,\ep)\ra\R$ are differentiable functions, satisfying 
\ea
{\d w_j\over\d t}&=a_j\,
\overline{w_1\cdots w_{j-1}w_{j+1}\cdots w_m},\quad j=1,\ldots,m,
\label{sy7eq4}\\
{\d u\over\d t}&=2\Re(w_1\cdots w_m),\quad\text{and}
\label{sy7eq5}\\
\ms{w_j}&=a_ju+1, \quad j=1,\ldots,m.
\label{sy7eq6}
\ea
If\/ \eq{sy7eq4} and\/ \eq{sy7eq5} hold for all\/ $t\in(-\ep,\ep)$ 
and\/ \eq{sy7eq6} holds for $t=0$, then \eq{sy7eq6} holds for all\/
$t$. Define a subset\/ $N$ of\/ $\C^m$ by
\e
\begin{split}
N=\Bigl\{\bigl(&r{\rm e}^{i\al_1}w_1(t),\ldots,r{\rm e}^{i\al_m}w_m(t)
\bigr):r>0,\quad t\in(-\ep,\ep),\\
&\al_j\in\R,\quad \al_1\!+\!\cdots\!+\!\al_m\!=\!0,\quad
a_1\al_1\!+\!\cdots\!+\!a_m\al_m\!=\!0\Bigr\}.
\end{split}
\label{sy7eq7}
\e
Then $N$ is a special Lagrangian submanifold in $\C^m$ with phase~$i^{m-2}$.
\label{sy7thm1}
\end{thm}

\subsection{Rewriting these equations}
\label{sy71}

We will eventually solve equations \eq{sy7eq4}--\eq{sy7eq6} of Theorem
\ref{sy7thm1} fairly explicitly. First, we find a way to rewrite them 
using fewer variables. Set
\begin{equation*}
w_j(t)={\rm e}^{i\th_j(t)}\sqrt{a_ju(t)+1} 
\end{equation*}
for differentiable functions $\th_1,\ldots,\th_m:(-\ep,\ep)\ra\R$. Define 
\begin{equation*}
\th=\th_1+\cdots+\th_m \quad\text{and}\quad Q(u)=\prod_{j=1}^m(a_ju+1).
\end{equation*}
Then calculation shows that \eq{sy7eq4}--\eq{sy7eq6} are equivalent to 
\ea
{\d u\over\d t}&=2\Re(w_1\cdots w_m)=2Q(u)^{1/2}\cos\th\quad\text{and}
\label{sy7eq8}\\
{\d\th_j\over\d t}&=-\,{a_jQ(u)^{1/2}\sin\th\over a_ju+1},
\quad j=1,\ldots,m.
\label{sy7eq9}
\ea
Summing \eq{sy7eq9} from $j=1$ to $m$ gives
\begin{equation*}
{\d\th\over\d t}=-\,Q(u)^{1/2}\sin\th\sum_{j=1}^m{a_j\over a_ju+1}.
\end{equation*}

Now although the point $(w_1,\ldots,w_m)$ is determined by the
$m+1$ real variables $\th_1,\ldots,\th_m$ and $u$, we are actually
only interested in the $\R_+\!\t G$-orbit ${\cal O}_t$ of $(w_1,\ldots,w_m)$.
It is not difficult to show that
\begin{align*}
{\cal O}_t=\Bigl\{\bigl(&r{\rm e}^{i\al_1}\sqrt{a_1u+1},\ldots,
r{\rm e}^{i\al_m}\sqrt{a_mu+1}\bigr):r>0,\\
&\al_j\in\R,\quad \al_1+\cdots+\al_m=\th,\quad
a_1\al_1+\cdots+a_m\al_m=\psi\Bigr\},
\end{align*}
where $\psi=a_1\th_1+\cdots+a_m\th_m$. Thus, ${\cal O}_t$ depends only
on the three real variables $\th,\psi$ and $u$. Furthermore,
\eq{sy7eq9} shows that $\psi$ evolves by
\begin{equation*}
{\d\psi\over\d t}=-\,Q(u)^{1/2}\sin\th\sum_{j=1}^m{a_j^2\over a_ju+1}.
\end{equation*}
Thus we may rewrite Theorem \ref{sy7thm1} as follows:

\begin{thm} Let\/ $u,\th,\psi:(-\ep,\ep)\ra\R$ be
differentiable functions satisfying 
\ea
{\d u\over\d t}&=2Q(u)^{1/2}\cos\th,
\label{sy7eq10}\\
{\d\th\over\d t}&=-\,Q(u)^{1/2}\sin\th\sum_{j=1}^m{a_j\over a_ju+1}
\label{sy7eq11}\\
\text{and}\quad
{\d\psi\over\d t}&=-\,Q(u)^{1/2}\sin\th\sum_{j=1}^m{a_j^2\over a_ju+1},
\label{sy7eq12}
\ea
such that\/ $a_ju(t)+1>0$ for $j=1,\ldots,m$ and\/ $t\in(-\ep,\ep)$.
Define a subset\/ $N$ of\/ $\C^m$ to be
\e
\begin{split}
\Bigl\{\bigl(&r{\rm e}^{i\al_1}\sqrt{a_1u(t)\!+\!1},\ldots,
r{\rm e}^{i\al_m}\sqrt{a_mu(t)\!+\!1}\,\bigr):r>0,\;\> 
t\in(-\ep,\ep),\\
&\al_j\in\R,\;\> \al_1+\cdots+\al_m=\th(t),\;\>
a_1\al_1+\cdots+a_m\al_m=\psi(t)\Bigr\}.
\end{split}
\label{sy7eq13}
\e
Then $N$ is a special Lagrangian submanifold in $\C^m$ with phase~$i^{m-2}$.
\label{sy7thm2}
\end{thm}

This is a significant simplification, as Theorem \ref{sy7thm1} was
written in terms of an o.d.e.\ in $m$ complex variables $w_1,\ldots,w_m$, 
but we have reduced this to only 3 real variables $u,\th,\psi$. In fact
we can show that $u$ and $\th$ are dependent in a simple way, and so
reduce the number of real variables to two.

Suppose for the moment that $\sin\th(0)\ne 0$, and divide \eq{sy7eq10} 
by \eq{sy7eq11}. This gives an expression for ${\d u\over\d\th}$, 
eliminating $t$. Separating variables shows that
\begin{equation*}
\int_{u(0)}^{u(t)}\sum_{j=1}^m{a_j\over a_ju+1}\,\d u=
-2\int_{\th(0)}^{\th(t)}\cot\th\,\d\th,
\end{equation*}
which integrates explicitly to
\begin{equation*}
\log Q(u)=-2\log\sin\th+C
\end{equation*}
for all $t\in(-\ep,\ep)$, for some $C\in\R$. So exponentiating 
gives~$Q(u)\sin^2\th\equiv {\rm e}^C>0$.  

If on the other hand $\sin\th(0)=0$ then \eq{sy7eq11} shows
that $\th$ is constant, so $Q(u)\sin^2\th\equiv 0$. In both
cases we see that $Q(u)\sin^2\th$ is constant, so its square 
root $Q(u)^{1/2}\sin\th$ is also constant, as it is continuous.
Thus we have $Q(u)^{1/2}\sin\th\equiv A$ for some~$A\in\R$.

This simplifies \eq{sy7eq11} and \eq{sy7eq12}, as we can replace the 
factor $Q(u)^{1/2}\sin\th$ by $A$. Also, from \eq{sy7eq10} we obtain
\begin{equation*}
\Bigl({\d u\over\d t}\Bigr)^2=
4Q(u)\cos^2\th=4\bigl(Q(u)-Q(u)\sin^2\th\bigr)=4\bigl(Q(u)-A^2\bigr).
\end{equation*}
Thus we have proved:

\begin{prop} In the situation of Theorem \ref{sy7thm2} we have
\e
Q(u)^{1/2}\sin\th\equiv A
\label{sy7eq14}
\e
for some $A\in\R$ and all\/ $t\in(-\ep,\ep)$, and\/ 
\eq{sy7eq10}--\eq{sy7eq12} are equivalent to
\begin{gather}
\Bigl({\d u\over\d t}\Bigr)^2=4\bigl(Q(u)-A^2\bigr),
\label{sy7eq15}\\
{\d\th\over\d t}=-A\sum_{j=1}^m{a_j\over a_ju+1}\quad\text{and}\quad
{\d\psi\over\d t}=-A\sum_{j=1}^m{a_j^2\over a_ju+1}.
\label{sy7eq16}
\end{gather}
\label{sy7prop1}
\end{prop}

The following lemma pins down the ranges of $u$ and~$A$.

\begin{lem} In the situation above, $u$ lies in $(-a_m^{-1},-a_1^{-1})$ 
for all\/ $t\in(-\ep,\ep)$, and\/ $0\le A^2\le Q(u)\le 1$. Also, $A$ 
lies in~$[-1,1]$.
\label{sy7lem1}
\end{lem}

\begin{proof} As by assumption $a_1,\ldots,a_m$ are not all zero with
$a_1\le\cdots\le a_m$ and $a_1+\cdots+a_m=0$, we see that $a_1<0$
and $a_m>0$, so the interval $(-a_m^{-1},-a_1^{-1})$ is well-defined. 
Since $a_ju+1>0$ for $j=1,\ldots,m$, we have $a_1u+1>0$ and $a_mu+1>0$.
Therefore $u$ is confined to~$(-a_m^{-1},-a_1^{-1})$.

Now $Q'(0)=\sum_{j=1}^ma_j=0$, so $Q$ has a turning point at 0. 
As the roots of $Q$ are $-a_1^{-1},\ldots,-a_m^{-1}$, they are 
all real, and there are none in $(-a_m^{-1},-a_1^{-1})$. So 0 is 
the only turning point of $Q$ in this interval. Thus $Q\le 1$ in 
$(-a_m^{-1},-a_1^{-1})$, as $Q(0)=1$. Equation \eq{sy7eq14} shows 
that $Q(u)\ge A^2\ge 0$ for $t\in(-\ep,\ep)$. As $Q(u)\le 1$, this 
gives $A^2\le 1$, so~$A\in[-1,1]$.
\end{proof}

As by changing the sign of one of the complex coordinates $z_j$ we
change the sign of $A$, we lose little by restricting our attention
to $A\in[0,1]$. This is equivalent to supposing that $\sin\th(0)\ge 0$.
To describe the SL $m$-folds $N$ of Theorem \ref{sy7thm2} in more detail, 
we shall divide into the three cases
\begin{itemize}
\item[(a)] $A=0$, 
\item[(b)] $A=1$, and
\item[(c)] $A\in(0,1)$.
\end{itemize}
We deal with each case separately.

\subsection{Case (a): $A=0$}
\label{sy72}

Suppose $A=0$ in Theorem \ref{sy7thm2} and Proposition \ref{sy7prop1}.
Then \eq{sy7eq16} shows that $\th$ and $\psi$ are constant, and
\eq{sy7eq14} that $\sin\th=0$, so $\th\equiv n\pi$ for some integer 
$n$, and \eq{sy7eq10} becomes ${\d u\over\d t}=2(-1)^nQ(u)^{1/2}$. So 
$u$ is monotone increasing or decreasing in $t$, and fills out some 
open interval in~$\R$.

The possible range of $u$ is $(-a_m^{-1},-a_1^{-1})$. Fixing 
$\th\equiv\psi\equiv 0$ and letting $u$ take its maximum range, 
Theorem \ref{sy7thm2} gives

\begin{prop} Let\/ $a_1\le\cdots\le a_m$ be integers, not all zero, 
with\/ $a_1+\cdots+a_m=0$. Define a subset\/ $N$ of\/ $\C^m$ to be
\e
\begin{split}
\Bigl\{\bigl(&r{\rm e}^{i\al_1}\sqrt{a_1u\!+\!1},\ldots,
r{\rm e}^{i\al_m}\sqrt{a_mu\!+\!1}\,\bigr):r>0,\quad
u\in(-a_m^{-1},-a_1^{-1}), \\
&\al_j\in\R,\;\> \al_1+\cdots+\al_m=0,\quad
a_1\al_1+\cdots+a_m\al_m=0\Bigr\}.
\end{split}
\label{sy7eq18}
\e
Then $N$ is a nonsingular SL submanifold in $\C^m$ 
with phase~$i^{m-2}$.
\label{sy7prop2}
\end{prop}

This special Lagrangian $m$-fold $N$ is not closed in $\C^m$. Its 
closure $\bar N$ is given by replacing the conditions $r>0$ and
$u\in(-a_m^{-1},-a_1^{-1})$ in \eq{sy7eq18} by $r\ge 0$ and 
$u\in[-a_m^{-1},-a_1^{-1}]$. Thus $\bar N\sm N$ is the disjoint 
union of 0 and two $\R_+\t G$-orbits ${\cal O}_+$, ${\cal O}_-$ 
corresponding to $u=-a_1^{-1}$ and~$u=-a_m^{-1}$.

Suppose ${\cal O}_\pm$ have dimension $m-1$. Then $\bar N$ is a 
{\it submanifold with boundary}, with an isolated singular point 
at 0. Furthermore, by Theorem \ref{sy6thm1} there must exist 
locally unique, cohomogeneity one special Lagrangian cones with 
phase $i^{m-2}$ containing ${\cal O}_+$ and ${\cal O}_-$. That is, 
we can extend $N$ as a special Lagrangian submanifold beyond the 
boundary components~${\cal O}_\pm$.

The condition that $\dim{\cal O}_+=m-1$ turns out to be that 
$a_1<a_2$, and the condition that $\dim{\cal O}_-=m-1$ turns out 
to be that $a_{m-1}<a_m$. It is then not difficult to prove the 
following result.

\begin{prop} Let\/ $a_1<a_2\le\cdots\le a_{m-1}<a_m$ be integers, 
not all zero, with\/ $a_1+\cdots+a_m=0$. Define a subset\/ $N$ of\/ 
$\C^m$ to be
\e
\begin{split}
\Bigl\{\bigl(&\pm r{\rm e}^{i\al_1}\sqrt{a_1u\!+\!1},
r{\rm e}^{i\al_2}\sqrt{a_2u\!+\!1},\ldots,
r{\rm e}^{i\al_{m-1}}\sqrt{a_{m-1}u\!+\!1},\\
&\pm r{\rm e}^{i\al_m}\sqrt{a_mu\!+\!1}\,\bigr):
r\ge 0,\;\> u\in[-a_m^{-1},-a_1^{-1}],\\
&\al_j\in\R,\;\> \al_1+\cdots+\al_m=0,\;\>
a_1\al_1+\cdots+a_m\al_m=0\Bigr\}.
\end{split}
\label{sy7eq19}
\e
Then $N$ is a closed, embedded SL cone in $\C^m$ 
with phase $i^{m-2}$, with an isolated singular point at\/~$0$.
\label{sy7prop3}
\end{prop}

The significant idea involved here is that when $u=-a_1^{-1}$,
the first complex coordinate $r{\rm e}^{i\al_1}\sqrt{a_1u\!+\!1}$
in \eq{sy7eq18} becomes zero. The natural way to extend $N$ beyond 
${\cal O}_+$ is to change the sign of the square root $\sqrt{a_1u\!+\!1}$.
Similarly, when $u=-a_m^{-1}$ the last complex coordinate becomes 
zero, and we can extend $N$ beyond ${\cal O}_-$ by changing
the sign of $\sqrt{a_mu\!+\!1}$. The $m$-fold $N$ of \eq{sy7eq19} 
is actually a closed loop of $\R^+\t G$-orbits, and so is a cone 
on a torus~$T^{m-1}$.

We have not yet considered the possibilities that ${\cal O}_+$ or 
${\cal O}_-$ have dimension less than $m-1$. Then we cannot use
Theorem \ref{sy6thm1} to extend $N$ past ${\cal O}_+$ or ${\cal O}_-$.
However, it turns out that if $\dim{\cal O}_+=m-2$, then $\bar N$ is 
already nonsingular without boundary at ${\cal O}_+$, and similarly
for ${\cal O}_-$. The local picture is of a family of circles collapsing 
to a point in $\R^2$. However, if $\dim{\cal O}_+\le m-3$ then $\bar N$
is unavoidably singular on~${\cal O}_+$.

The condition that $\dim{\cal O}_+=m-2$ is that $a_1=a_2<a_3$, and 
the condition that $\dim{\cal O}_-=m-2$ is that $a_{m-2}<a_{m-1}=a_m$.
When $\dim{\cal O}_+=m-1$ and $\dim{\cal O}_-=m-2$, we can prove the 
following result by the same method as Proposition~\ref{sy7prop3}.

\begin{prop} Let\/ $a_1<a_2\le\cdots\le a_{m-2}<a_{m-1}=a_m$ be integers, 
not all zero, with\/ $a_1+\cdots+a_m=0$. Define a subset\/ $N$ of\/ 
$\C^m$ to be
\begin{align*}
\Bigl\{\bigl(&\pm r{\rm e}^{i\al_1}\sqrt{a_1u\!+\!1},
r{\rm e}^{i\al_2}\sqrt{a_2u\!+\!1},\ldots,
r{\rm e}^{i\al_m}\sqrt{a_mu\!+\!1}\,\bigr):u\in[-a_m^{-1},-a_1^{-1}],\\
&r\ge 0,\;\> \al_j\in\R,\;\> \al_1+\cdots+\al_m=0,\;\>
a_1\al_1+\cdots+a_m\al_m=0\Bigr\}.
\end{align*}
Then $N$ is a closed, embedded SL cone in $\C^m$ 
with phase $i^{m-2}$, with an isolated singular point at\/~$0$.
\label{sy7prop4}
\end{prop}

In this case $N$ is a closed interval of $\R_+\t G$-orbits rather 
than a loop, and is topologically a cone on ${\cal S}^2\t T^{m-3}$.
When $m=3$ we must take $a_1=-2$ and $a_2=a_3=1$, and we just get a 
copy of $\R^3$ in $\C^3$. But for $m\ge 4$, the proposition gives 
new, interesting SL $m$-folds in~$\C^m$.

Similarly, when $\dim{\cal O}_+=\dim{\cal O}_-=m-2$, we prove:

\begin{prop} Let\/ $a_1=a_2<a_3\le\cdots\le a_{m-2}<a_{m-1}=a_m$ be 
integers, not all zero, with\/ $a_1+\cdots+a_m=0$. Define a subset\/ 
$N$ of\/ $\C^m$ to be
\begin{align*}
\Bigl\{\bigl(&r{\rm e}^{i\al_1}\sqrt{a_1u\!+\!1},\ldots,
r{\rm e}^{i\al_m}\sqrt{a_mu\!+\!1}\,\bigr):r\ge 0,\;\>
u\in[-a_m^{-1},-a_1^{-1}],\\
&\al_j\in\R,\;\> \al_1+\cdots+\al_m=0,\;\>
a_1\al_1+\cdots+a_m\al_m=0\Bigr\}.
\end{align*}
Then $N$ is a closed, embedded SL cone in $\C^m$ 
with phase $i^{m-2}$, with an isolated singular point at\/~$0$.
\label{sy7prop5}
\end{prop}

Again, $N$ is a closed interval of $\R_+\t G$-orbits rather than 
a loop, and is topologically a cone on ${\cal S}^3\t T^{m-4}$.
The conditions on $a_1,\ldots,a_m$ do not hold when $m=3$. When
$m=4$ we must take $a_1=a_2=-1$ and $a_3=a_4=1$, and we just get
a copy of $\R^4$ in $\C^4$. But for $m\ge 5$, the proposition
yields new, interesting SL $m$-folds in~$\C^m$.

All these examples can also be constructed using the `perpendicular
symmetry' construction to be described in \S\ref{sy9}. Specifically,
in Proposition \ref{sy9prop} we set $n=m$ and $c=0$, and define $G$ 
as above. The proposition then yields a special Lagrangian
submanifold $N_0$ in $\C^m$, of which the $m$-folds $N$ constructed
above will be subsets.

\subsection{Case (b): $A=1$}
\label{sy73}

Next suppose that $A=1$ in Theorem \ref{sy7thm2} and Proposition
\ref{sy7prop1}. As $A^2\le Q(u)\le 1$ this gives $Q(u)\equiv 1$,
which forces $u\equiv 0$, as $Q(0)=1$ is the strict maximum of
$Q$ in the permitted interval $(-a_m^{-1},-a_1^{-1})$. Then
\eq{sy7eq14} gives $\sin\th\equiv 1$, so that $\th\equiv\pi/2$. 
The second equation of \eq{sy7eq16} then becomes ${\d\psi\over\d t}=
-\sum_{j=1}^ma_j^2$. Hence solutions to \eq{sy7eq10}--\eq{sy7eq12} 
exist for all $t\in\R$, and are given by
\begin{equation*}
u(t)=0,\quad
\th(t)={\pi\over 2} \quad\text{and}\quad
\psi(t)=\psi(0)-t\sum_{j=1}^ma_j^2,
\end{equation*}
for $\psi(0)\in\R$.

Now $\psi$ takes all values in $\R$. Thus the equation
$a_1\al_1+\cdots+a_m\al_m=\psi(t)$ in \eq{sy7eq13} is
actually no restriction. So the SL $m$-fold $N$ of 
\eq{sy7eq13} reduces to
\e
N=\Bigl\{\bigl(r{\rm e}^{i\al_1},\ldots,
r{\rm e}^{i\al_m}\bigr):r>0,\;\>
\al_j\in\R,\;\> \al_1+\cdots+\al_m={\pi\over 2}\Bigr\}.
\label{sy7eq20}
\e
This $N$ is entirely independent of $a_1,\ldots,a_m$. It has
generalized symmetry group $\GSym(N)=\R_+\t\U(1)^{m-1}$, which 
acts transitively, and symmetry group $\Sym(N)=\U(1)^{m-1}$, 
which acts with cohomogeneity one. 

Now we have already studied SL submanifolds of $\C^m$ on which 
$\U(1)^{m-1}$ acts with cohomogeneity one in Example \ref{cmex1}. 
In fact $N$ is half of the SL $m$-fold $L_{0,\ldots,0}$ of Example 
\ref{cmex1}, rotated to give it phase $i^{m-2}$ rather than~1.

\subsection{Case (c): $A\in(0,1)$, local treatment}
\label{sy74}

We shall discuss case (c) above from two points of view. Firstly,
when $\cos\th>0$ for $t\in(-\ep,\ep)$ we will find a more explicit 
expression for the manifold $N$ of Theorem \ref{sy7thm2} by 
eliminating $t$, and writing $\th$ and $\psi$ as functions of $u$. 
Then in \S\ref{sy75} we will discuss the behaviour of equations 
\eq{sy7eq10}--\eq{sy7eq12} for $t\in\R$ rather than $(-\ep,\ep)$, 
and show that they admit periodic solutions.

We would like to write the SL $m$-fold $N$ of Theorem \ref{sy7thm2} 
in as simple and explicit a way as possible. One way of doing this 
is to eliminate $t$, and write everything instead as a function of 
the variable $u$. Now ${\d u\over\d t}$ has the same sign as $\cos\th$ 
by \eq{sy7eq10}. Thus, if $\cos\th$ changes sign in $(-\ep,\ep)$ then 
we cannot write $t$ as a function of $u$, but if $\cos\th$ has constant 
sign then we can.

Let us {\it assume} that $\cos\th>0$ for for all $t\in(-\ep,\ep)$. 
Then \eq{sy7eq15} gives ${\d u\over\d t}=2\sqrt{Q(u)-A^2}>0$, and
integrating gives
\e
\int_{u(0)}^{u(t)}{\d u\over 2\sqrt{Q(u)-A^2}}=\int_0^t\d t=t.
\label{sy7eq21}
\e
This equation defines $u$ implicitly as a function of $t$.
The integral on the left hand side is called an {\it elliptic 
integral} if $m=3$ or 4, and a {\it hyperelliptic integral}
if~$m\ge 5$.

From \eq{sy7eq16} and the expression above for 
${\d u\over\d t}$ we get
\begin{equation*}
{\d\psi\over\d u}=-\,{A\over 2\sqrt{Q(u)-A^2}}
\sum_{j=1}^m{a_j^2\over a_ju+1}.
\end{equation*}
Integrating this gives an expression for $\psi$ in terms of $u$.
Setting $u_0=u(0)$ and $u_{\pm\ep}=u(\pm\ep)$, we have proved:

\begin{thm} Under the assumptions above, the SL $m$-fold\/ $N$ of 
Theorem \ref{sy7thm2} may be written as
\begin{align*}
\Bigl\{\bigl(&r{\rm e}^{i\al_1}\sqrt{a_1u\!+\!1},\ldots,
r{\rm e}^{i\al_m}\sqrt{a_mu\!+\!1}\,\bigr):r>0,\;\> 
u\in\bigl(u_{-\ep},u_\ep\bigr),\\
&\al_j\in\R,\;\> \al_1+\cdots+\al_m=\th(u),\;\>
a_1\al_1+\cdots+a_m\al_m=\psi(u)\Bigr\},
\end{align*}
where $\th(u)$ and\/ $\psi(u)$ are given by
$\th(u)=\sin^{-1}\bigl(AQ(u)^{-1/2}\bigr)$ and
\e
\psi(u)=\psi(u_0)-{A\over 2}\int_{u_0}^u{\d v\over\sqrt{Q(v)-A^2}}
\sum_{j=1}^m{a_j^2\over a_jv+1}.
\label{sy7eq22}
\e
\label{sy7thm3}
\end{thm}

This is a reasonably explicit expression for $N$. The integral
defining $\psi$ probably cannot be simplified any further without 
making special assumptions about~$a_1,\ldots,a_m$.

\subsection{Case (c): $A\in(0,1)$, global behaviour}
\label{sy75}

Next we study the global behaviour of solutions to equations 
\eq{sy7eq10}--\eq{sy7eq12} of Theorem \ref{sy7thm2} when $A\in(0,1)$.
We begin with a preliminary lemma on the range of~$u$.

\begin{lem} Suppose $A\in(0,1)$. Then there exist constants $\al<0,\be>0$ 
and\/ $\ga>0$ with\/ $Q(\al)=Q(\be)=A^2$ and\/ $Q(u)>A^2$ for 
$u\in(\al,\be)$, such that for all\/ $t$ for which solutions $u,\th$ 
exist to \eq{sy7eq10}--\eq{sy7eq11}, we have $u(t)\in[\al,\be]$ and\/ 
$a_ju(t)+1\ge\ga>0$ for~$j=1,\ldots,m$.
\label{sy7lem2}
\end{lem}

\begin{proof} From the proof of Lemma \ref{sy7lem1} we know that $Q$ 
strictly increases from 0 to 1 in $[-a_m^{-1},0]$, and so as $A\in(0,1)$ 
there is a unique $\al\in(-a_m^{-1},0)$ with $Q(\al)=A^2$. Similarly, 
there is a unique $\be\in(0,-a_1^{-1})$ with $Q(\be)=A^2$. Clearly 
$Q(u)>A^2$ for~$u\in(\al,\be)$. 

Also, if $u$ lies in $(-a_m^{-1},-a_1^{-1})$ then $Q(u)\ge A^2$ 
if and only if $u\in[\al,\be]$. But we know from Lemma \ref{sy7lem1}
that if $u,\th$ are solutions to \eq{sy7eq10}--\eq{sy7eq11} then $u$ 
is confined to $(-a_m^{-1},-a_1^{-1})$, and $Q(u)\ge A^2$. Thus 
$u(t)\in[\al,\be]$ for all $t$ for which the solutions exist, as we 
have to prove.

As $u$ is confined to $[\al,\be]$, there exists $K>0$ with 
$a_ju+1\le K$ for $j=1,\ldots,m$ and all $t$ for which the 
solution exists. Thus
\begin{equation*}
0<A^2\le Q(u)=(a_ju+1)\prod_{\substack{1\le i\le m\\ i\ne j}}(a_iu+1)
\le (a_ju+1)K^{m-1},
\end{equation*}
so that $a_ju(t)+1\ge\ga=A^2K^{1-m}>0$ for~$j=1,\ldots,m$.
\end{proof}

Now we can show that solutions exist for all $t\in\R$, and $u,\th$
are periodic.

\begin{prop} Let\/ $u(0)$, $\th(0)$ and\/ $\psi(0)$ be given, with\/ 
$a_ju(0)+1>0$ for $j=1,\ldots,m$ and\/ $A=Q(u(0))^{1/2}\sin\th(0)\in(0,1)$.
Then there exist unique solutions $u(t)$, $\th(t)$ and\/ $\psi(t)$ to 
equations \eq{sy7eq10}--\eq{sy7eq12} of Theorem \ref{sy7thm2} for all\/ 
$t\in\R$, with these values at\/ $t=0$. Furthermore $u$ and\/ $\th$ are 
nonconstant and periodic with period\/ $T>0$, and there exists $\Psi>0$ 
with\/ $\psi(t+T)=\psi(t)-\Psi$ for all\/~$t\in\R$.
\label{sy7prop7}
\end{prop}

\begin{proof} The only way for solutions of \eq{sy7eq10}--\eq{sy7eq12} 
to become singular is for some $a_ju+1$ to become zero, or for 
$u\ra\pm\iy$. As neither of these can happen by Lemma \ref{sy7lem2}, 
the solutions must exist for all $t\in\R$. Uniqueness of the solutions, 
with the given initial data, follows from standard results in 
differential equations.

Since $A>0$ we have $\sin\th>0$ by \eq{sy7eq14}. By \eq{sy7eq10},
${\d u\over\d t}=0$ if and only if $\cos\th=0$, that is, if and
only if $\sin\th=1$. By \eq{sy7eq14} this happens exactly when 
$Q(u)=A^2$. But $u$ is confined to $[\al,\be]$ by Lemma 
\ref{sy7lem2}, where $Q(\al)=Q(\be)=A^2$ and $Q(u)>A^2$ for 
$u\in(\al,\be)$. Hence ${\d u\over\d t}=0$ if and only if 
$u=\al$ or~$u=\be$.

So we see that $u$ must cycle up and down between $\al$ and $\be$,
turning only at $\al$ and $\be$. Using the ideas of \S\ref{sy74},
we see from \eq{sy7eq21} that the time taken for $u$ to increase
from $\al$ to $\be$ is 
\begin{equation*}
\ha T=\int_\al^\be{\d u\over 2\sqrt{Q(u)-A^2}},
\end{equation*}
which is finite, as $Q(u)-A^2$ has only single roots at $u=\al,\be$.
Similarly, the time for $u$ to decrease from $\be$ to $\al$ is
also $\ha T$. Thus, the time taken for $u$ to start at $\al$,
increase to $\be$, and decrease back to $\al$, is $T$. That is, 
$u$ undergoes periodic oscillations with period~$T$.

Equation \eq{sy7eq16} shows that ${\d\psi\over\d t}$ is also
periodic with period $T$, and also that ${\d\psi\over\d t}<0$,
as $A>0$. It is then easy to see that $\psi(t+T)=\psi(t)-\Psi$
for all $t$ and $\Psi>0$ given by~$\Psi=\psi(0)-\psi(T)$.
\end{proof}

Next we show that if $\Psi$ is a rational multiple $a/b$ of $2\pi$ 
then the family of $\R_+\!\t G$-orbits ${\cal O}_t$ is periodic, with 
period $bT$, and $N$ is a cone on~$T^{m-1}$.

\begin{prop} In the situation of Proposition \ref{sy7prop7},
suppose that\/ $\Psi=2\pi q$ for $q\in\Q$. Define a subset\/ 
$N$ of\/ $\C^m$ to be
\e
\begin{split}
\Bigl\{\bigl(&r{\rm e}^{i\al_1}\sqrt{a_1u(t)\!+\!1},\ldots,
r{\rm e}^{i\al_m}\sqrt{a_mu(t)\!+\!1}\,\bigr):r\ge 0,\;\> 
t\in\R,\\
&\al_j\in\R,\;\> \al_1+\cdots+\al_m=\th(t),\;\>
a_1\al_1+\cdots+a_m\al_m=\psi(t)\Bigr\}.
\end{split}
\label{sy7eq23}
\e
Then $N$ is a closed, embedded SL cone in $\C^m$ with phase $i^{m-2}$, 
which is topologically a cone on $T^{m-1}$, and has just one singular 
point at\/~$0$.
\label{sy7prop8}
\end{prop}

\begin{proof} Let $q=a/b$, for $a,b$ positive integers with $\hcf(a,b)=1$.
Then as $u(t+T)=u(t)$, $\th(t+T)=\th(t)$ and $\psi(t+T)=\psi(t)-\Psi$, we 
have $u(t+bT)=u(t)$, $\th(t+bT)=\th(t)$ and $\psi(t+bT)=\psi(t)-2\pi a$.
Now the definition \eq{sy7eq13} of the SL $m$-fold
$N$ constructed in Theorem \ref{sy7thm2} actually depends only on 
${\rm e}^{i\th}$ and ${\rm e}^{i\psi}$ rather than on $\th$ and $\psi$, 
as $a_1,\ldots,a_m$ are integers. Thus, adding an integer multiple of 
$2\pi$ to $\psi$ makes no difference to~$N$.

This means that the the 1-parameter family of $\R_+\!\t G$-orbits
${\cal O}_t$ which make up $N$ satisfies ${\cal O}_{t+bT}={\cal O}_t$ 
for all $t$, and is periodic, with period $bT$. The definition of $N$ 
given in \eq{sy7eq23} differs from \eq{sy7eq13} in that $r\ge 0$ rather 
than $r>0$, and $t\in\R$ rather than $t\in(-\ep,\ep)$. We allow $r=0$ 
to add the point 0 to $N$, which makes $N$ {\it closed}. We allow 
$t\in\R$ because $u,\th$ and $\psi$ exist for all $t\in\R$ by
Proposition \ref{sy7prop7}. In fact, as $t$ has period $bT$, we
would get the same $N$ if we replaced $t\in\R$ by~$t\in[0,bT)$.

So $N$ is a closed loop of $\R_+\!\t G$-orbits ${\cal O}_t$, together
with the point zero. It is clear that $N$ is closed and, at least 
as an immersed submanifold, it is topologically a cone on $T^{m-1}$, 
with just one singular point at 0. So we need only show that $N$ is 
embedded. If two $\R_+\!\t G$-orbits ${\cal O}_t,{\cal O}_{t'}$ intersect, then
they are the same. But from Theorem \ref{sy6thm1}, an $\R_+\!\t G$-orbit
locally determines $N$ uniquely. Thus, our loop of $\R_+\!\t G$-orbits 
cannot cross itself, and must be embedded.
\end{proof}

The proposition is potentially interesting because closed special 
Lagrangian cones in $\C^m$ with isolated singular points are natural
models for singularities of SL $m$-folds in Calabi--Yau
$m$-folds. We will now investigate the range of $\Psi$, and hence
show that the proposition yields very many such cones.

By the reasoning used to prove \eq{sy7eq22}, we can show that
\e
\Psi=\int_\al^\be{\d v\over\sqrt{A^{-2}Q(v)-1}}
\sum_{j=1}^m{a_j^2\over a_jv+1}.
\label{sy7eq24}
\e
As Lemma \ref{sy7lem2} defines $\al,\be$ in terms of $Q(u)$
and $A$, we see that $\Psi$ depends only on $a_1,\ldots,a_m$ and 
$A$, and not on the initial data $u(0),\th(0)$ and~$\psi(0)$.

Up to now we have regarded $A$ as a function of $u(0),\th(0)$ and
$\psi(0)$. We now change our point of view. Lemma \ref{sy7lem2}
defines $\al,\be$ depending on $A\in(0,1)$. Given any $A\in(0,1)$, 
set $u(0)=\al$, $\th(0)=\pi/2$ and $\psi(0)=0$. This is a valid set 
of initial data, and yields this value of $A$. Thus $A$ can take any 
value in $(0,1)$, and we can regard $u(0),\th(0),\psi(0)$ and $\Psi$ 
as functions of $A$. We now calculate the limit of $\Psi$ 
as $A$ approaches 0 or~1.

\begin{prop} In the situation above, $\Psi:(0,1)\ra(0,\iy)$ is a real 
analytic function of\/ $A$, and satisfies $\Psi(A)\ra\pi(a_m-a_1)$ as 
$A\ra 0$ and $\Psi(A)\ra\pi\bigl(2\sum_{j=1}^ma_j^2\bigr)^{1/2}$ as~$A\ra 1$.
\label{sy7prop9}
\end{prop}

\begin{proof} It is obvious from \eq{sy7eq24} and the definition
of $\al$ and $\be$ that $\Psi$ is real analytic. As $A\ra 0$, we 
have $\al\ra -1/a_m$ and $\be\ra -1/a_1$. Also, as $A\ra 0$ the 
factor $(A^{-2}Q(v)-1)^{-1/2}$ in \eq{sy7eq24} tends to zero, 
except near $\al$ and $\be$. Hence, as $A\ra 0$, the integrand 
in \eq{sy7eq24} gets large near $\al\approx -1/a_m$ and 
$\be\approx -1/a_1$, and very close to zero in between. 

So to understand $\Psi$ as $A\ra 0$, it is enough to study the integral 
\eq{sy7eq24} near $\al$ and $\be$. We shall model it at $\al$. Suppose 
$a_{m-k}<a_{m-k+1}=\cdots=a_m$, so that $a_m$ has multiplicity $k$. 
Then near $v=-1/a_m$ we have
\begin{equation*}
Q(v)\approx C(v+a_m^{-1})^k,
\qquad\text{where}\qquad
C=a_m^k\prod_{j=1}^{m-k}(1-a_j/a_m).
\end{equation*}
Since $A^2=Q(\al)$ this gives $A^2\approx C(\al+a_m^{-1})^k$, so 
that~$\al\approx A^{2/k}C^{-1/k}-a_m^{-1}$.

Therefore, when $v\approx\al$ we have
\begin{equation*}
A^{-2}Q(v)-1\approx A^{-2}C(v+a_m^{-1})^k-1
\quad\text{and}\quad
\sum_{j=1}^m{a_j^2\over a_jv+1}\approx{ka_m\over v+a_m^{-1}},
\end{equation*}
taking only the highest-order terms. Thus, when $A$ is small we
see that
\begin{align*}
&\int_\al^0{\d v\over\sqrt{A^{-2}Q(v)-1}}
\sum_{j=1}^m{a_j^2\over a_jv+1}\approx\\
&\int_{A^{2/k}C^{-1/k}-a_m^{-1}}^0{\d v\over\sqrt{A^{-2}C(v+a_m^{-1})^k-1}}
\cdot{ka_m\over v+a_m^{-1}}.
\end{align*}
Changing variables to $w=\sqrt{A^{-2}C(v+a_m^{-1})^k-1}$, after
some surprising cancellations we get
\begin{equation*}
\int_{A^{2/k}C^{-1/k}-a_m^{-1}}^0{\d v\over\sqrt{A^{-2}C(v+a_m^{-1})^k-1}}
\cdot{ka_m\over v+a_m^{-1}}
\approx\int_0^\iy{2a_m\d w\over w^2+1}=\pi a_m,
\end{equation*}
where we have approximated the second integral by replacing the
upper limit $\sqrt{A^{-2}Ca_m^{-k}-1}$ by $\iy$. Hence, for small 
$A$ we have 
\begin{equation*}
\int_\al^0{\d v\over\sqrt{A^{-2}Q(v)-1}}
\sum_{j=1}^m{a_j^2\over a_jv+1}\approx \pi a_m,
\end{equation*}
and similarly
\begin{equation*}
\int_0^\be{\d v\over\sqrt{A^{-2}Q(v)-1}}
\sum_{j=1}^m{a_j^2\over a_jv+1}\approx -\pi a_1,
\end{equation*}
so that $\Psi(A)\ra\pi(a_m-a_1)$ as~$A\ra 0$.

Next consider the behaviour of $\Psi$ as $A\ra 1$. When $A$ 
is close to 1, $u$ is small and $\th$ is close to $\pi/2$. So 
write $\th={\pi\over 2}+\phi$, for $\phi$ small. Then, 
setting~$Q(u)\approx 1$,
\begin{equation*}
\cos\th\approx -\phi,\;\>
\sin\th\approx 1,\;\>
\sum_{j=1}^m{a_j\over a_ju+1}\approx -u\sum_{j=1}^ma_j^2
\;\>\text{and}\;\>
\sum_{j=1}^m{a_j^2\over a_ju+1}\approx \sum_{j=1}^ma_j^2,
\end{equation*}
taking only the highest order terms, equations \eq{sy7eq10}--\eq{sy7eq12} 
become
\begin{equation*}
{\d u\over\d t}\approx -2\phi,\quad
{\d\phi\over\d t}\approx u\sum_{j=1}^ma_j^2,
\quad\text{and}\quad
{\d\psi\over\d t}\approx -\sum_{j=1}^ma_j^2.
\end{equation*}
The first two of these equations show that $u$ and $\phi$ 
undergo approximately simple harmonic oscillations with period 
$T=2\pi\bigl(2\sum_{j=1}^ma_j^2\bigr)^{-1/2}$. Then the third
equation shows that $\Psi\approx -{\d\psi\over\d t}T$, as
${\d\psi\over\d t}$ is approximately constant, and so
$\Psi\ra\pi\bigl(2\sum_{j=1}^ma_j^2\bigr)^{1/2}$ as~$A\ra 1$.
\end{proof}

Now the limits $\pi(a_m-a_1)$ and $\pi\bigl(2\sum_{j=1}^ma_j^2\bigr)^{1/2}$
of $\Psi(A)$ as $A\ra 0,1$ satisfy
\begin{equation*}
\Bigl[\pi\bigl(2{\ts\sum}_{j=1}^ma_j^2\bigr)^{1/2}\Bigr]^2-
\bigl[\pi(a_m-a_1)\bigr]^2=
\pi^2(a_1+a_m)^2+2\pi^2(a_2^2+\cdots+a_{m-1}^2).
\end{equation*}
Thus $\pi(a_m-a_1)\le\pi(2\sum_{j=1}^ma_j^2)^{1/2}$, with equality 
if and only if $a_1+a_m=0$ and $a_2=\cdots=a_{m-1}=0$. As the $a_j$ 
are integers with highest common factor 1, and $a_1\le\cdots\le a_m$,
these conditions imply that $a_1=-1$, $a_2=\cdots=a_{m-1}=0$ and~$a_m=1$.

Therefore we have two cases:
\begin{itemize}
\item[(i)] $(a_1,\ldots,a_m)\ne (-1,0,\ldots,0,1)$, and
$\lim_{A\ra 0}\Psi(A)<\lim_{A\ra 1}\Psi(A)$, or
\item[(ii)] $(a_1,\ldots,a_m)=(-1,0,\ldots,0,1)$, and
$\lim_{A\ra 0}\Psi(A)=\lim_{A\ra 1}\Psi(A)=2\pi$.
\end{itemize}
In case (i), as $\lim_{A\ra 0}\Psi(A)<\lim_{A\ra 1}\Psi(A)$ we
see that $\Psi$ is not constant, and as it is real analytic it 
can have only finitely many stationary points in $(0,1)$. So we deduce:

\begin{cor} Suppose $(a_1,\ldots,a_m)\ne (-1,0,\ldots,0,1)$. Then for 
a countable dense subset of\/ $A\in(0,1)$ we have~$\Psi(A)\in 2\pi\Q$.
\label{sy7cor}
\end{cor}

In case (ii), we can solve equations \eq{sy7eq4}--\eq{sy7eq6} 
completely. For as $a_2=\cdots=a_{m-1}=0$ we see from \eq{sy7eq4} 
that $w_2,\ldots,w_{m-1}$ are constant, and from \eq{sy7eq6} that 
$\md{w_2}=\cdots=\md{w_{m-1}}=1$. Applying a diagonal matrix in 
$\SU(m)$, we may choose $w_2=\cdots=w_{m-1}=1$. Then \eq{sy7eq4} 
and \eq{sy7eq6} reduce to
\begin{equation*}
{\d w_1\over\d t}=-\bar w_m,\quad
{\d w_m\over\d t}=\bar w_1\quad\text{and}\quad
\ms{w_1}+\ms{w_m}=2,
\end{equation*}
which have solutions
\begin{equation*}
w_1=B{\rm e}^{it}+C{\rm e}^{-it}\quad\text{and}\quad
w_m=i\bar B{\rm e}^{-it}-i\bar C{\rm e}^{it},
\end{equation*}
for $B,C\in\C$ with $\ms{B}+\ms{C}=1$. It is easy to show that
$A=2\Im(\bar BC)\in[-1,1]$, and that $\Psi(A)=2\pi$ for all~$A$.

The special Lagrangian $m$-fold $N$ of \eq{sy7eq7} is thus
\begin{align*}
N=\Bigl\{&\bigl(r{\rm e}^{i\al_1}(B{\rm e}^{it}+C{\rm e}^{-it}),
r{\rm e}^{i\al_2},\ldots,r{\rm e}^{i\al_{m-1}},
r{\rm e}^{i\al_m}(i\bar B{\rm e}^{-it}-i\bar C{\rm e}^{it})\bigr):\\
&r>0,\;\> t\in\R,\;\> \al_j\in\R,\;\> \al_1+\cdots+\al_m=0,\;\> 
\al_1=\al_m\Bigr\}.
\end{align*}
Now this is the result of applying the $\SU(m)$ transformation
\begin{equation*}
(z_1,\ldots,z_m)\mapsto
(Bz_1-iCz_m,z_2,\ldots,z_{m-1},-i\bar Cz_1+\bar Bz_m)
\end{equation*}
of $\C^m$ to the special Lagrangian cone
\begin{align*}
N'=\Bigl\{&\bigl(r{\rm e}^{i(\al_1+t)},
r{\rm e}^{i\al_2},\ldots,r{\rm e}^{i\al_{m-1}},
ir{\rm e}^{i(\al_m-t)}\bigr):\\
&r>0,\quad t\in\R,\quad \al_j\in\R,\quad \al_1+\cdots+\al_m=0,\quad
\al_1=\al_m\Bigr\}.
\end{align*}
But this is identical to the cone defined in \eq{sy7eq20}, which 
we saw in \S\ref{sy73} has transitive generalized symmetry group 
$\GSym(N)=\R_+\t\U(1)^{m-1}$, and symmetry group~$\Sym(N)=\U(1)^{m-1}$.

Here is how to interpret this. Since $\Psi(A)=2\pi$ for all $A$, 
we would expect Proposition \ref{sy7prop8} to yield a 1-parameter 
family of distinct $\U(1)^{m-2}$-invariant SL cones on $T^{m-1}$ 
in $\C^m$, parametrized by $A\in(0,1)$, and with the same symmetry 
group $\U(1)^{m-2}$. But in fact these SL cones are all isomorphic
under transformations in $\SU(m)$, and have symmetry group 
$\U(1)^{m-1}$ rather than~$\U(1)^{m-2}$.

One consequence of this is that there are {\it no} SL cones $N$ on 
$T^{m-1}$ in $\C^m$ with $\Sym^0(N)$ equal to this particular symmetry 
group $G=\U(1)^{m-2}$ in $\SU(m)$, since any SL cone symmetric under 
this group $\U(1)^{m-2}$ is also symmetric under a larger group
$\U(1)^{m-1}$ in~$\SU(m)$.

Drawing together much of the work above, in particular Proposition 
\ref{sy7prop8} and Corollary \ref{sy7cor}, we have the main result
of this section.

\begin{thm} Let\/ $G$ be a Lie subgroup of\/ $\SU(m)$ isomorphic to 
$\U(1)^{m-2}$, for some $m\ge 3$. If\/ $G$ is conjugate in $\SU(m)$
to the group
\e
\begin{split}
\Bigl\{&(z_1,\ldots,z_m)\mapsto
({\rm e}^{i\al_1}z_1,\ldots,{\rm e}^{i\al_m}z_m):\\
&\al_j\in\R,\quad \al_1+\cdots+\al_m=0,\quad \al_1=\al_m\Bigr\},
\end{split}
\label{sy7eq25}
\e
then every $G$-invariant SL cone in $\C^m$ locally has symmetry
group~$\U(1)^{m-1}$.

Otherwise there exists a countably infinite family of distinct, 
closed, embedded SL cones $N$ in $\C^m$ with\/ $\Sym^0(N)=G$,
each of which is topologically a cone on $T^{m-1}$, with just 
one singular point at\/~$0$.
\label{sy7thm4}
\end{thm}

Here is what we mean by saying the construction yields a countable 
family of distinct cones. Fixing a subgroup $G\cong\U(1)^{m-2}$ in
$\SU(m)$, we consider two $G$-invariant SL cones $N_1,N_2$ to be 
equivalent if there exists a $G$-equivariant isometric isomorphism 
of $\C^m$ identifying $N_1$ and $N_2$. The theorem says that for 
$G$ not conjugate to \eq{sy7eq25} there are countably many distinct 
equivalence classes of $G$-equivariant SL cones.

More generally, we could ask about the classification of the 
SL cones constructed in the second part of the theorem for all 
$\U(1)^{m-2}$-subgroups $G$, up to isometric isomorphisms of 
$\C^m$. Let $N_1,N_2$ be two such cones with groups $G_1,G_2$, 
and suppose $N_1$ and $N_2$ are identified by an isometric 
isomorphism of $\C^m$. Then $G_1$ and $G_2$ are conjugate in 
$\SU(m)$, since $G_j=\Sym^0(N_j)$. 

Thus, for each of the countable number of conjugacy classes of 
$\U(1)^{m-2}$-subgroups $G$ in $\SU(m)$ not conjugate to \eq{sy7eq25}, 
we have found a countable family of distinct $G$-invariant special
Lagrangian $T^{m-1}$-cones in $\C^m$, up to isometric isomorphisms 
of $\C^m$. One reason these cones are interesting is as local models 
for singularities of SL $m$-folds in Calabi--Yau $m$-folds. 

To see how big the whole family is, here is a crude `parameter count'. 
Up to isomorphism, $N$ depends on integers $a_1\le\cdots\le a_m$ 
with $a_1+\cdots+a_m=0$ and $\hcf(a_1,\ldots,a_m)=1$ and $A\in(0,1)$ 
with $\Psi(A)=2\pi q$, for $q\in\Q$. The most obvious thing to do is 
to set $q=a/b$ with $\hcf(a,b)=1$, and say that $N$ depends on the 
$m+1$ integers $a_1,\ldots,a_{m-1},a$ and~$b$.

However, this is probably not the best point of view. Instead,
we should drop the condition $\hcf(a_1,\ldots,a_m)=1$, and
replace the $a_j$ by $\ti a_j=a_jb$ for $j=1,\ldots,m$. With
these new values we get $\Psi(A)=2\pi a$, so that $\Psi$ lies 
in $2\pi\Z$ rather than $2\pi\Q$, and we can say that $N$ depends 
on the $m$ integers $\ti a_1,\ldots,\ti a_{m-1},a$, which have
highest common factor~1.

A partial version of the case $m=3$ of Theorem \ref{sy7thm4} was
first due to Mark Haskins \cite[Th.~C]{Hask}. Haskins does not 
study the periodicity conditions in the case $A\in(0,1)$, but only 
when $A=0$. The author's treatment was completed independently 
somewhat later, and uses different methods.

\subsection{Relation to integrable systems}
\label{sy76}

We saw above that equation \eq{sy7eq4} has very nice behaviour ---
various quantities are conserved, and one can say a lot about
the solutions, even writing them explicitly using elliptic
integrals. The reason for this is that \eq{sy7eq4} is {\it 
completely integrable}, as we will now show.

An introduction to integrable systems is given in Hitchin, 
Segal and Ward \cite{HSW}. From \cite[Def.~6.1, p.~49]{HSW}
and \cite[Def.~2.7, p.~61]{HSW}, we may define a {\it completely
integrable Hamiltonian system} to be a $2m$-dimensional symplectic 
manifold $(M,\om)$ with a Hamiltonian $H:M\ra\R$, such that there
exist $m$ conserved quantities $p_1,\ldots,p_m:H\ra\R$ which are
independent almost everywhere and satisfy
\e
\{H,p_i\}\equiv \{p_i,p_j\}\equiv 0\quad\text{for $i,j=1,\ldots,m$},
\label{sy7eq26}
\e
where $\{\,,\,\}$ is the Poisson bracket on~$M$.

Choose $a_1,\ldots,a_m\in\R\sm\{0\}$. Let $M$ be $\C^m$ with 
coordinates $(w_1,\ldots,w_m)$, and define a symplectic form 
$\om$ and a Hamiltonian function $H$ on $M$ by
\begin{equation*}
\om={i\over 2}\sum_{j=1}^ma_j^{-1}\d w_j\w\d\bar w_j
\quad\text{and}\quad
H:(w_1,\ldots,w_m)\mapsto 2\Im(w_1\cdots w_m).
\end{equation*}

The equations of motion of this Hamiltonian system are
easily shown to be
\begin{equation*}
{\d w_j\over\d t}=a_j\,
\overline{w_1\cdots w_{j-1}w_{j+1}\cdots w_m},\quad j=1,\ldots,m,
\end{equation*}
as in \eq{sy7eq4}. Define functions $p_1,\ldots,p_m:M\ra\R$ by
\begin{equation*}
p_j=a_m\ms{w_j}-a_j\ms{w_m} \quad\text{for $j=1,\ldots,m-1$, and}\quad
p_m=\Im(w_1\cdots w_m).
\end{equation*}
It is easy to show that $p_1,\ldots,p_m$ are conserved, independent 
almost everywhere, and satisfy \eq{sy7eq26}. Thus \eq{sy7eq4} is 
indeed a completely integrable Hamiltonian system.

This proof assumes that $a_1,\ldots,a_m$ are nonzero, but not 
that $a_1+\cdots+a_m=0$, which is the case in \eq{sy7eq4}. If
some $a_j$ is zero then ${\d w_j\over\d t}=0$, so $w_j$ is
constant. Fixing all the $w_j$ with $a_j=0$, we can regard 
\eq{sy7eq4} as an o.d.e.\ in the remaining $w_k$, which is 
then integrable as above.

\section{$\U(1)$-invariant SL cones in $\C^3$}
\label{sy8}

We now specialize to the case $m=3$ in the situation of
\S\ref{sy7}, so that we are studying SL cones $N$ in $\C^3$ 
invariant under a $\U(1)$ subgroup of $\SU(3)$. In this
case we can improve the treatment of \S\ref{sy7} in several
ways. In particular, we will define the group $G\cong\U(1)$ 
in a neater way, we will write down conformal coordinates on 
$N\cap{\cal S}^5$, and we will solve equations 
\eq{sy7eq10}--\eq{sy7eq12} explicitly in terms of Jacobi 
elliptic functions.

The material of this section has already been studied by
several authors, from different points of view. Probably the 
first were Castro and Urbano \cite{CaUr1}, who constructed
examples of minimal Lagrangian tori in $\CP^2$ using integrable
systems methods. Special Lagrangian $T^2$-cones in $\C^3$ can 
be reconstructed from Castro and Urbano's results. 

Later, Haskins \cite[\S 3--\S 5]{Hask} studied $\U(1)$-invariant 
special Lagrangian cones in $\C^3$, and most of this section 
overlaps with his work. In particular, the author learnt the 
material of \S\ref{sy82} from his papers.

From Theorem \ref{sy7thm1} we get:

\begin{thm} Let\/ $a_1,a_2,a_3\in\R$ be not all zero, with\/ 
$a_1+a_2+a_3=0$. Suppose $w_1,w_2,w_3:\R\ra\C$ and\/ $u:\R\ra\R$ satisfy
\ea
{\d w_1\over\d t}&=a_1\,\overline{w_2w_3},\quad
{\d w_2\over\d t}=a_2\,\overline{w_3w_1},\quad
{\d w_3\over\d t}=a_3\,\overline{w_1w_2},
\label{sy8eq9}\\
{\d u\over\d t}&=2\Re(w_1w_2w_3),\quad\text{and}
\label{sy8eq10}\\
\ms{w_j}&=a_ju+1\quad\text{for $j=1,2,3$.}
\label{sy8eq11}
\ea
If\/ \eq{sy8eq9} and\/ \eq{sy8eq10} hold for all\/ $t$ and\/
\eq{sy8eq11} holds for $t=0$, then \eq{sy8eq11} holds for all\/ 
$t$. Define a subset\/ $N$ of\/ $\C^3$ by
\e
\begin{split}
N=\Bigl\{\bigl(&r{\rm e}^{i\al_1}w_1(t),r{\rm e}^{i\al_2}w_2(t),
r{\rm e}^{i\al_3}w_3(t)\bigr):r>0,\quad t\in\R,\\
&\al_j\in\R,\quad \al_1+\al_2+\al_3=0,\quad
a_1\al_1+a_2\al_2+a_3\al_3=0\Bigr\}.
\end{split}
\label{sy8eq12}
\e
Then $N$ is a special Lagrangian $3$-fold in $\C^3$ with phase~$i$.
\label{sy8thm1}
\end{thm}

Here we have replaced the domain $(-\ep,\ep)$ of $w_j$ and $u$ by $\R$, 
as the results of \S\ref{sy7} imply that solutions of \eq{sy8eq9} exist 
for all $t\in\R$. We will find it convenient later to rewrite the theorem 
in terms of $b_1,b_2,b_3$, where
\begin{equation*}
b_1={\ts\frac{1}{\sqrt{3}}}(a_2-a_3),\quad 
b_2={\ts\frac{1}{\sqrt{3}}}(a_3-a_1)\quad\text{and}\quad
b_3={\ts\frac{1}{\sqrt{3}}}(a_1-a_2).
\end{equation*}
Then $b_1+b_2+b_3=0$, and as $a_1+a_2+a_3=0$ we find that
\e
a_1={\ts\frac{1}{\sqrt{3}}}(b_3-b_2),\quad 
a_2={\ts\frac{1}{\sqrt{3}}}(b_1-b_3) \quad\text{and}\quad 
a_3={\ts\frac{1}{\sqrt{3}}}(b_2-b_1).
\label{sy8eq13}
\e
It is easy to show that $\al_1,\al_2,\al_3\in\R$ satisfy
$\al_1+\al_2+\al_3=0$ and $a_1\al_1+a_2\al_2+a_3\al_3=0$ if and
only if $\al_j=b_js$ for some $s\in\R$. Also, \eq{sy8eq11} is
equivalent to $\ms{w_1}+\ms{w_2}+\ms{w_3}=3$ and $b_1\ms{w_1}
+b_2\ms{w_2}+b_3\ms{w_3}=0$. Thus, Theorem \ref{sy8thm1} becomes:

\begin{thm} Let\/ $b_1,b_2,b_3\in\R$ be not all zero, with\/ 
$b_1+b_2+b_3=0$. Define $a_1,a_2,a_3$ by \eq{sy8eq13}. Suppose 
$w_1,w_2,w_3:\R\ra\C\sm\{0\}$ satisfy
\begin{gather}
{\d w_1\over\d t}=a_1\overline{w_2w_3},\quad
{\d w_2\over\d t}=a_2\overline{w_3w_1}\quad\text{and}\quad 
{\d w_3\over\d t}=a_3\overline{w_1w_2},
\label{sy8eq14}\\
\ms{w_1}+\ms{w_2}+\ms{w_3}=3 \quad\text{and}\quad
b_1\ms{w_1}+b_2\ms{w_2}+b_3\ms{w_3}=0.
\label{sy8eq15}
\end{gather}
If\/ \eq{sy8eq14} holds for all\/ $t$ and\/ \eq{sy8eq15} holds
for $t=0$, then \eq{sy8eq15} holds for all\/ $t$. Define a subset\/ 
$N$ of\/ $\C^3$ by
\e
N=\Bigl\{\bigl(r{\rm e}^{ib_1s}w_1(t),\ldots,
r{\rm e}^{ib_3s}w_3(t)\bigr):r>0,\;\> 
s\in\R,\;\> t\in\R\Bigr\}.
\label{sy8eq16}
\e
Then $N$ is a special Lagrangian $3$-fold in $\C^3$ with phase~$i$.
\label{sy8thm2}
\end{thm}

We shall show that $(s,t)$ are {\it conformal coordinates} on the 
unit sphere $N\cap{\cal S}^5$ in $N$, and that the corresponding map 
$\Phi:\R^3\ra{\cal S}^5$ is {\it harmonic}.

\begin{prop} In the situation of Theorem \ref{sy8thm2}, define 
$\Phi:\R^2\ra{\cal S}^5$ by 
\e
\Phi:(s,t)\mapsto{\ts\frac{1}{\sqrt{3}}}\bigl({\rm e}^{ib_1s}w_1(t),
{\rm e}^{ib_2s}w_2(t),{\rm e}^{ib_3s}w_3(t)\bigr),
\label{sy8eq17}
\e
where ${\cal S}^5$ is the unit sphere in $\C^3$. Then $\Phi$ is a 
conformal harmonic map.
\label{sy8prop1}
\end{prop}

\begin{proof} From \eq{sy8eq14} and \eq{sy8eq17} we see that
\begin{align*}
{\pd\Phi\over\pd s}&={\ts\frac{1}{\sqrt{3}}}\bigl(ib_1{\rm e}^{ib_1s}w_1,
ib_2{\rm e}^{ib_2s}w_2,ib_3{\rm e}^{ib_3s}w_3\bigr)\quad\text{and}\\
{\pd\Phi\over\pd t}&={\ts\frac{1}{\sqrt{3}}}\bigl(a_1{\rm e}^{ib_1s}
\overline{w_2w_3},a_2{\rm e}^{ib_2s}\overline{w_3w_1},
a_3{\rm e}^{ib_3s}\overline{w_1w_2}\bigr).
\end{align*}
Now $\Phi$ is conformal if and only if ${\pd\Phi\over\pd s}$ and
${\pd\Phi\over\pd t}$ are orthogonal and the same length. But
\begin{equation*}
g\bigl({\ts{\pd\Phi\over\pd s}},{\ts{\pd\Phi\over\pd t}}
\bigr)=-{\ts\frac{1}{3}}\bigl[a_1b_1+a_2b_2+a_3b_3\bigr]\Im(w_1w_2w_3)=0,
\end{equation*}
as $a_1b_1+a_2b_2+a_3b_3=0$, so they are orthogonal. Also
\begin{align*}
\bms{\ts{\pd\Phi\over\pd s}}&={\ts\frac{1}{3}}
\bigl(b_1^2\ms{w_1}+b_2^2\ms{w_2}+b_3^2\ms{w_3}\bigr)\quad\text{and}\\
\bms{\ts{\pd\Phi\over\pd t}}&={\ts\frac{1}{3}}\bigl(a_1^2
\ms{w_2}\ms{w_3}+a_2^2\ms{w_3}\ms{w_1}+a_3^2\ms{w_1}\ms{w_2}\bigr).
\end{align*}
One can then prove from equations \eq{sy8eq13}, \eq{sy8eq15} and
$b_1+b_2+b_3=0$ that $\ms{{\pd\Phi\over\pd s}}=\ms{{\pd\Phi\over\pd t}}$, 
and thus $\Phi$ is conformal.

Now $N$ is a cone in $\C^3$, and is minimal because any calibrated
submanifold is minimal. Therefore the intersection $N\cap{\cal S}^5$
is minimal in ${\cal S}^5$. But $N\cap{\cal S}^5$ is the image of
$\Phi$. Hence $\Phi:\R^2\ra{\cal S}^5$ is a conformal map with
minimal image. But it is well known in the field of harmonic maps
that a conformal map from a Riemann surface is harmonic if and
only if it has minimal image. Thus $\Phi$ is harmonic.
\end{proof}

Using the method of \S\ref{sy71} to rewrite the $w_j$ in terms of
real variables $u,\th_j$, as in Theorem \ref{sy7thm2} and Proposition
\ref{sy7prop1} we get

\begin{thm} In the situation of Theorem \ref{sy8thm2} the functions 
$w_1,w_2,w_3$ may be written $w_j={\rm e}^{i\th_j}\sqrt{a_ju+1}$, for 
$u,\th_1,\th_2,\th_3:\R\ra\R$. Define
\begin{gather*}
Q(u)=(a_1u+1)(a_2u+1)(a_3u+1),\\
\th=\th_1+\th_2+\th_3\qquad\text{and}\qquad
\psi=a_1\th_1+a_2\th_2+a_2\th_3.
\end{gather*}
Then $Q(u)^{1/2}\sin\th\equiv A$ for some $A\in\R$, and\/
$u,\th_j,\th$ and\/ $\psi$ satisfy
\e
\begin{aligned}
\Bigl({\d u\over\d t}\Bigr)^2&=4\bigl(Q(u)-A^2\bigr),\quad&
{\d\th_j\over\d t}&=-\,{a_jA\over a_ju+1},\\
{\d\th\over\d t}&=-A\sum_{j=1}^3{a_j\over a_ju+1},\quad&
{\d\psi\over\d t}&=-A\sum_{j=1}^3{a_j^2\over a_ju+1}.
\end{aligned}
\label{sy8eq18}
\e
The phase $i$ special Lagrangian $3$-fold\/ $N$ of\/ \eq{sy8eq16} 
may also be written
\begin{align*}
N=\Bigl\{\bigl(&r{\rm e}^{i\al_1}\sqrt{a_1u(t)+1},
\ldots,r{\rm e}^{i\al_3}\sqrt{a_3u(t)+1}\bigr):
r>0,\;\> t\in\R\\
&\al_j\in\R,\;\> \al_1+\al_2+\al_3=\th(t),\;\>
a_1\al_1+a_2\al_2+a_3\al_3=\psi(t)\Bigr\}.
\end{align*}
\label{sy8thm3}
\end{thm}

Our next theorem follows from the case $m=3$ of Theorem~\ref{sy7thm4}.

\begin{thm} Let\/ $b_1,b_2,b_3$ be distinct integers with highest
common factor $1$ and\/ $b_1+b_2+b_3=0$, and define a subgroup
$G\cong\U(1)$ of\/ $\SU(3)$ to be the set of transformations 
\begin{equation*}
(z_1,z_2,z_3)\mapsto({\rm e}^{ib_1s}z_1,{\rm e}^{ib_2s}z_2,
{\rm e}^{ib_3s}z_3)
\end{equation*}
for $s\in\R$. Then there exists a countably infinite family of 
distinct, closed, embedded special Lagrangian cones $N$ in $\C^3$ 
with\/ $\Sym^0(N)=G$, each of which is topologically a cone on 
$T^2$, and has just one singular point at\/~$0$.
\label{sy8thm4}
\end{thm}

Here we have tidied the theorem up by defining the group $G$ using
$b_1,b_2,b_3$ instead of $a_1,a_2,a_3$, and requiring the $b_j$ to
be integers rather than the $a_j$. The condition that $G$ should 
not be conjugate in $\SU(3)$ to the group \eq{sy7eq25} turns
out to be that $b_1,b_2,b_3$ are distinct. 

The sense in which the cones are distinct was explained after
Theorem \ref{sy7thm4}. Two sets of integers $b_1,b_2,b_3$ and
$b_1',b_2',b_3'$ produce isomorphic families of SL cones in $\C^3$ if 
and only if the corresponding groups $G,G'$ are conjugate in $\SU(3)$. 
This happens if and only if $b_j'=b_{\sigma(j)}$ for $j=1,2,3$ or 
$b_j'=-b_{\sigma(j)}$ for $j=1,2,3$, for some~$\sigma\in S_3$.

\subsection{The Jacobi elliptic functions}
\label{sy81}

We now give a brief introduction to the {\it Jacobi elliptic functions},
which we will use in \S\ref{sy82} to solve \eq{sy8eq18} explicitly. The 
following material can be found in Chandrasekharan \cite[Ch.~VII]{Chan}. 
For each $k\in[0,1]$, the Jacobi elliptic functions $\sn(t,k)$, $\cn(t,k)$, 
$\dn(t,k)$ with modulus $k$ are the unique solutions to the o.d.e.s
\begin{align}
\bigl({\ts{\d\over\d t}}\sn(t,k)\bigr)^2&=\bigl(1-\sn^2(t,k)\bigr)
\bigl(1-k^2\sn^2(t,k)\bigr),
\label{jacobieq1}\\
\bigl({\ts{\d\over\d t}}\cn(t,k)\bigr)^2&=\bigl(1-\cn^2(t,k)\bigr)
\bigl(1-k^2+k^2\cn^2(t,k)\bigr),
\label{jacobieq2}\\
\bigl({\ts{\d\over\d t}}\dn(t,k)\bigr)^2&=-\bigl(1-\dn^2(t,k)\bigr)
\bigl(1-k^2-\dn^2(t,k)\bigr),
\label{jacobieq3}
\end{align}
with initial conditions
\e
\begin{aligned}
\sn(0,k)&=0,\quad & \cn(0,k)&=1,\quad & \dn(0,k)&=1,\\
{\ts{\d\over\d t}}\sn(0,k)&=1,\quad&
{\ts{\d\over\d t}}\cn(0,k)&=0,\quad&
{\ts{\d\over\d t}}\dn(0,k)&=0.
\end{aligned}
\label{jacobieq4}
\e
They satisfy the identities
\e
\sn^2(t,k)+\cn^2(t,k)=1 \quad\text{and}\quad k^2\sn^2(t,k)+\dn^2(t,k)=1,
\label{jacobieq5}
\e
and the differential equations
\e
\begin{gathered}
{\ts{\d\over\d t}}\sn(t,k)=\cn(t,k)\dn(t,k),\qquad
{\ts{\d\over\d t}}\cn(t,k)=-\sn(t,k)\dn(t,k)\\
\text{and}\qquad {\ts{\d\over\d t}}\dn(t,k)=-k^2\sn(t,k)\cn(t,k).
\end{gathered}
\label{jacobieq6}
\e
When $k=0$ or 1 they reduce to trigonometric functions:
\e
\begin{aligned}
\sn(t,0)&=\sin t,&\quad \cn(t,0)&=\cos t,\quad \dn(t,0)=1,\\
\sn(t,1)&=\tanh t&,\quad \cn(t,1)&=\dn(t,1)=\sech t.
\end{aligned}
\label{jacobieq7}
\e
For $k\in[0,1)$ the Jacobi elliptic functions are
periodic in $t$, with $\sn(t,k)$ and $\cn(t,k)$ of
period $4K(k)$ and $\dn(t,k)$ of period $2K(k)$, where
\e
K(k)=\int_0^{\pi/2}{\d x\over \sqrt{1-k^2\sin^2x}}\,.
\label{jacobieq8}
\e

\subsection{Explicit solution using Jacobi elliptic functions}
\label{sy82}

Following Haskins \cite[\S 4]{Hask}, we shall solve \eq{sy8eq18} 
fairly explicitly. The answer depends on the order of $a_1,a_2,a_3$
and 0. For simplicity, we shall suppose that $a_2\le a_1<0<a_3$.
The solutions for the other possible orders may be obtained by
permutations of 1,2,3 and sign changes, in an obvious way.

Now $Q(u)-A^2$ has three real roots $\ga_1,\ga_2,\ga_3$, which 
may be ordered so that $\ga_1\ge\ga_2\ge 0\ge\ga_3$. Then 
$Q(u)-A^2=a_1a_2a_3(u-\ga_1)(u-\ga_2)(u-\ga_3)$, and the first 
equation of \eq{sy8eq18} becomes
\begin{equation*}
\Bigl({\d u\over\d t}\Bigr)^2=4a_1a_2a_3(u-\ga_1)(u-\ga_2)(u-\ga_3).
\end{equation*}
As in \cite[Prop.~4.2]{Hask} the solutions are $u(t)=\ga_3+
(\ga_2-\ga_3)\sn^2(at+c,b)$, where $a^2=a_1a_2a_3(\ga_1-\ga_3)$, 
$b^2=(\ga_2-\ga_3)/(\ga_1-\ga_3)$, $c\in\R$, and $\sn(\,,\,)$ is the 
Jacobi sn-noidal function. This can easily be verified using~\eq{jacobieq1}.

Substituting into the other three equations of \eq{sy8eq18}
gives explicit expressions for ${\d\th_j\over\d t}$,
${\d\th\over\d t}$ and ${\d\psi\over\d t}$, so we obtain
$\th_j,\th$ and $\psi$ by integration. We have proved:

\begin{prop} Under the assumptions above, the solutions of\/ 
\eq{sy8eq18} are
\ea
u(t)&=\ga_3+(\ga_2-\ga_3)\sn^2(at+c,b)\quad\text{and}
\label{sy8eq19}\\
\th_j(t)&=\th_j(0)-A
\int_0^t{\d t\over\ga_3+{\ts{1\over a_j}}
+(\ga_2-\ga_3)\sn^2(at+c,b)},
\label{sy8eq20}
\ea
where $a^2=a_1a_2a_3(\ga_1-\ga_3)$, $b^2=(\ga_2-\ga_3)/(\ga_1-\ga_3)$
and\/~$c\in\R$.
\label{sy8prop2}
\end{prop}

When $A=0$ we can write $w_1,w_2,w_3$ explicitly in terms of $\sn,\cn$ 
and $\dn$. By \eq{sy8eq20} we may take $\th_j\equiv 0$, so that 
$w_1,w_2,w_3$ are real. The roots $\ga_1,\ga_2,\ga_3$ of $Q(u)-A^2$ are 
$-1/a_1,-1/a_2,-1/a_3$, and the orders $a_2\le a_1<0<a_3$ and 
$\ga_1\ge\ga_2\ge 0\ge\ga_3$ imply that $\ga_j=-1/a_j$. Thus, combining 
\eq{sy8eq11} with \eq{sy8eq19} gives explicit expressions for $w_j^2$, 
which by \eq{jacobieq5} and the definition of $b^2$ above reduce to
\begin{equation*}
w_1^2={\ts{a_3-a_1\over a_3}}\,\dn^2(at,b),\quad
w_2^2={\ts{a_3-a_2\over a_3}}\,\cn^2(at,b),\quad
w_3^2={\ts{a_2-a_3\over a_2}}\,\sn^2(at,b),
\end{equation*}
putting $c=0$ for simplicity. Hence, from Theorem \ref{sy8thm2} 
we deduce:

\begin{thm} Let\/ $b_1,b_2,b_3\in\Z$ satisfy $b_2>b_3>0>b_1$ and\/
$b_1+b_2+b_3=0$. Define $a_1,a_2,a_3$ by \eq{sy8eq13}, and\/
$a>0$ and\/ $b\in(0,1)$ by 
\e
a^2=a_2(a_1-a_3) \quad\text{and}\quad 
b^2={a_1(a_2-a_3)\over a_2(a_1-a_3)}.
\label{sy8eq21}
\e
Define a subset\/ $N$ of\/ $\C^3$ by
\e
\begin{split}
N=\Bigl\{
\bigl(&r{\rm e}^{ib_1s}
\bigl({\ts{a_3-a_1\over a_3}}\bigr)^{1/2}\dn(at,b),
r{\rm e}^{ib_2s}
\bigl({\ts{a_3-a_2\over a_3}}\bigr)^{1/2}\cn(at,b),\\
&r{\rm e}^{ib_3s}
\bigl({\ts{a_2-a_3\over a_2}}\bigr)^{1/2}\sn(at,b)
\bigr):r>0,\quad s,t\in\R\Bigr\}.
\end{split}
\label{sy8eq22}
\e
Then $N$ is a special Lagrangian cone on $T^2$ in $\C^3$ with 
phase $i$. Furthermore, $(s,t)$ are conformal coordinates 
on~$N\cap{\cal S}^5$.
\label{sy8thm5}
\end{thm}

Here the expressions for $a^2$ and $b^2$ come from Proposition
\ref{sy8prop2} by putting $\ga_j=-1/a_j$, and $b_2>b_3\ge 0>b_1$ 
is equivalent to the condition $a_2\le a_1<0<a_3$ above. The 
additional assumption $b_3>0$ ensures that $b\in(0,1)$, as $b=1$ 
if and only if $b_3=0$. We have taken $b_1,b_2,b_3\in\Z$ to make 
the $s$ coordinate periodic, with period $2\pi$. We know from 
\S\ref{sy81} that $\sn(at,b)$, $\cn(at,b)$ and $\dn(at,b)$ are 
periodic in $t$, as $b\in(0,1)$. Thus $N$ is indeed a cone on~$T^2$.

\subsection{Relation with harmonic tori in $\CP^2$ and ${\cal S}^5$}
\label{sy83}

In Proposition \ref{sy8prop1} we showed that each of the SL 3-folds 
$N$ in $\C^3$ constructed in Theorem \ref{sy8thm2} is the cone on the 
image of a {\it conformal harmonic map}\/ $\Phi:\R^2\ra{\cal S}^5$ 
defined by \eq{sy8eq17}. Then in \S\ref{sy82} we showed that the 
functions $w_1,w_2,w_3$ in \eq{sy8eq17} may be written explicitly 
in terms of the Jacobi elliptic functions.

Thus we have constructed a family of explicit conformal harmonic 
maps $\Phi:\R^2\ra{\cal S}^5$. Furthermore, as $N$ is Lagrangian, one 
can show that if $\pi:{\cal S}^5\ra\CP^2$ is the Hopf projection then
$\pi\circ\Phi$ is conformal and harmonic, so we also have a family of 
explicit conformal harmonic maps~$\Psi:\R^2\ra\CP^2$.

Now harmonic maps from Riemann surfaces into spheres and projective
spaces are an {\it integrable system}, and have been intensively
studied in the integrable systems literature. For an introduction
to the subject, see Fordy and Wood \cite{FoWo}, in particular the
articles by Bolton and Woodward \cite[p.~59--82]{FoWo}, McIntosh
\cite[p.~205--220]{FoWo} and Burstall and Pedit~\cite[p.~221--272]{FoWo}.

Therefore our examples can be analyzed from the integrable systems
point of view. We postpone this analysis to the sequel \cite{Joyc6}. 
In \cite[\S 5]{Joyc6} we shall realize the SL 3-folds of Theorem 
\ref{sy8thm2} as special cases of a more general construction of 
special Lagrangian cones in $\C^3$, which involves two commuting 
o.d.e.s. Generic cones arising from this construction have only 
discrete symmetry groups.

Then in \cite[\S 6]{Joyc6} we work through the integrable systems
framework for the corresponding family of harmonic maps $\Psi:\R^2
\ra\CP^2$, showing that they are generically superconformal of 
finite type, and determining their harmonic sequences, Toda 
solutions, algebras of polynomial Killing fields, and spectral 
curves. From the integrable systems point of view, Theorem 
\ref{sy8thm4} is interesting because it constructs a large 
family of {\it superconformal harmonic tori} in~$\CP^2$.

\section{Construction by `perpendicular symmetries'}
\label{sy9}

Now we explain a new construction of special Lagrangian submanifolds
in $\C^m$ beginning with a special Lagrangian $m$-fold $L$ with
`perpendicular symmetries', that is, vector fields in $\su(m)\lt\C^m$
which are perpendicular to $L$ at every point.

\begin{thm} Let\/ $G$ be a $k$-dimensional abelian Lie subgroup 
of\/ $\SU(m)\lt\C^m$ with Lie algebra $\g$, acting on $\C^m$ 
with moment map $\mu$, and let\/ $\phi:\g\ra C^\iy(T\C^m)$ be the 
corresponding action of\/ $\g$ on $\C^m$ by vector fields. Suppose 
$L$ is an SL submanifold of\/ $\C^m$, such that\/ 
$\phi(x)$ is normal to $L$ at\/ $L$ for every $x$ in $\g$. For 
each\/ $c\in\g^*$, define $N_c=G\cdot\bigl(L\cap\mu^{-1}(c)\bigr)$.
Then $N_c$ is special Lagrangian in $\C^m$, with phase $1$ if\/ $k$ 
is even, and phase $i$ if\/ $k$ is odd.
\label{sy9thm1}
\end{thm}

\begin{proof} We shall show that for each nonsingular point
$z$ in $N_c$, the tangent plane $T_zN_c$ is special Lagrangian.
But as $N_c$ is $G$-invariant and $G\subset\SU(m)\lt\C^m$, it 
is enough to verify this for one point in each orbit of $G$ in 
$N_c$. Thus we can restrict our attention to $z\in L\cap\mu^{-1}(c)$.
Then the condition for $z$ to be a nonsingular point of $N_c$ is that 
$z$ is a nonsingular point of $L$, and the vector fields $\phi(\g)$ 
are linearly independent at~$z$. 

Suppose these conditions hold. Choose a basis $x_1,\ldots,x_k$
of $\g$ such that $\phi(x_1)_z,\ldots,\phi(x_k)_z$ are orthonormal, 
which is possible by linear independence of $\phi(\g)$ at $z$. 
Now $\phi(x_j)_z$ is normal to $T_zL$, which is a special 
Lagrangian plane in $\C^m$. Therefore $I\bigl(\phi(x_j)_z\bigr)$ 
lies in $T_zL$, where $I$ is the complex structure on $\C^m$.
Hence $I\bigl(\phi(x_1)_z\bigr),\ldots,I\bigl(\phi(x_k)_z\bigr)$
are orthonormal in $T_zL$. Extend them to an orthonormal basis of
$T_zL$ with vectors $v_1,\ldots,v_{m-k}$, so that 
\e
T_zL=\Big\langle I\bigl(\phi(x_1)_z\bigr),\ldots,
I\bigl(\phi(x_k)_z\bigr),v_1,\ldots,v_{m-k}\Big\rangle.
\label{tzleq}
\e

Now the level sets of the moment map $\mu$ of $G$ are orthogonal 
to $I\bigl(\phi(x_j)_z\bigr)$ for all $j$. Hence 
$T_z\bigl(L\cap\mu^{-1}(c)\bigr)$ is the subspace of $T_zL$
orthogonal to $I\bigl(\phi(x_j)_z\bigr)$ for $j=1,\ldots,k$. Thus 
$T_z\bigl(L\cap\mu^{-1}(c)\bigr)=\an{v_1,\ldots,v_{m-k}}$. But 
$N_c=G\cdot\bigl(L\cap\mu^{-1}(c)\bigr)$, and so $T_zN_c$ is the 
span of $T_z\bigl(L\cap\mu^{-1}(c)\bigr)$ and $\phi(\g)_z$. Hence
\e
T_zN_c=\ban{\phi(x_1)_z,\ldots,\phi(x_k)_z,v_1,\ldots,v_{m-k}}.
\label{tznceq}
\e

Comparing \eq{tzleq} and \eq{tznceq} and remembering that the bases 
are orthonormal, we see that in effect we have orthogonal direct sums 
$T_zL=\R^k\op\R^{m-k}$ and $T_zN_c=I(\R^k)\op\R^{m-k}$. It is easy
to see that as $T_zL$ is an SL plane with phase 1,
this implies that $T_zN_c$ is an SL plane with phase
$i^k$ or $-i^k$, depending on the orientation chosen for $N_c$. Thus,
if $k$ is even then $N_c$ is special Lagrangian with phase 1, and
if $k$ is odd then $N_c$ is special Lagrangian with phase $i$, with
the appropriate orientation.
\end{proof}

The assumption that $G$ is abelian was not actually used in 
the proof; but it is implied by the hypotheses, which is why 
we put it in. In the situation of the theorem, suppose $\g$ 
is not abelian, and let $x,y\in\g$. Then $I(\phi(x))$,
$I(\phi(y))$ are vector fields on $\C^m$ tangent to $L$ at 
$L$. Hence the Lie bracket $\bigl[I(\phi(x)),I(\phi(y))\bigr]$ 
is also tangent to $L$ at $L$. But 
$\phi(x),\phi(y)$ are holomorphic, and 
so~$\bigl[I(\phi(x)),I(\phi(y))\bigr]=
-\bigl[\phi(x),\phi(y)\bigr]=-\phi\bigl([x,y]\bigr)$.

Therefore $\phi\bigl([x,y]\bigr)$ is tangent to $L$, but 
it is also perpendicular to $L$. So $\phi\bigl([x,y]\bigr)=0$ on 
$L$. This forces $\phi\bigl([x,y]\bigr)=0$ on $\C^m$, since 
otherwise $L$ lies in some affine $\C^{m'}\subset\C^m$ for $m'<m$, 
which contradicts $L$ being Lagrangian. If $\phi$ is effective
then we have shown that $[x,y]=0$ for all $x,y\in\g$, so that 
$\g$ is abelian, and thus $G$ is abelian as it is connected.

We now characterize the possibilities for $G$ and $L$ in 
Theorem~\ref{sy9thm1}.

\begin{thm} Let\/ $G$ be a connected abelian Lie subgroup of\/ 
$\SU(m)\lt\C^m$ acting on the affine space $\C^m$, let\/ $\g$
be the Lie algebra of\/ $G$, and let\/ $\phi:\g\ra C^\iy(T\C^m)$ 
be the corresponding action of\/ $\g$ on $\C^m$ by vector fields. 

Then there exists an SL submanifold\/ $L$ of\/ $\C^m$ 
such that\/ $\phi(x)$ is normal to $L$ at\/ $L$ for every $x$ in $\g$, 
if and only if the following conditions hold:
\begin{itemize}
\setlength{\itemsep}{0pt}
\setlength{\parsep}{0pt}
\item[{\rm(i)}] There exists a $G$-invariant affine isomorphism
$\C^m\cong\C^{a_1}\t\cdots\t\C^{a_{n+2}}$;
\item[{\rm(ii)}] $L$ is a subset of the product manifold\/ 
$L_1\t L_2\t \cdots\t L_{n+2}$, where $L_j$ is a special
Lagrangian submanifold of\/~$\C^{a_j}$;
\item[{\rm(iii)}] For $j=1,\ldots,n$, $L_j$ is a cone in
$\C^{a_j}$, and each\/ $\ga\in G$ acts on $\C^{a_j}$ by 
multiplication by ${\rm e}^{i\th_j}$; 
\item[{\rm(iv)}] $L_{n+1}=\R^{a_{n+1}}$ in $\C^{a_{n+1}}$, and\/ 
$G$ acts on $\C^{a_{n+1}}$ by translations in the direction of\/ 
$I(\R^{a_{n+1}})$; and
\item[{\rm(v)}] $G$ acts trivially on $\C^{a_{n+2}}$.
\end{itemize}
\label{sy9thm2}
\end{thm}

\begin{proof} For simplicity, we first suppose that $G$ lies in the 
subgroup $\SU(m)$ of $\SU(m)\lt\C^m$, and treat $\C^m$ as a vector
space rather than an affine space. Then $\g$ is an abelian Lie 
subalgebra of $\su(m)$, which we may regard as a vector space of 
commuting matrices. By standard results in linear algebra, we may 
decompose the complex vector space $\C^m$ into a direct sum of 
{\it eigenspaces} of the action of~$\g$. 

Actually, one usually considers the eigenspaces of a single matrix, 
rather than of a vector space of commuting matrices. But the 
eigenspace decomposition of $\C^m$ under a generic element of $\g$ 
is the same as its decomposition under $\g$, so the two points of 
view are equivalent.

Let us write the eigenspace decomposition as 
\e
\C^m=\C^{a_1}\op\cdots\op\C^{a_n}\op V, 
\label{sy9eq1}
\e
where $\C^{a_1},\ldots,\C^{a_n}$ are nonzero eigenspaces of $\g$ in 
$\C^m$ with distinct, nonzero eigenvalues in $i\g^*$, and $V$ is the 
zero eigenspace of $\g$. The decomposition \eq{sy9eq1} is unique up 
to the order of the subspaces $\C^{a_1},\ldots,\C^{a_n}$, and is 
orthogonal as~$\g\subset{\mathfrak u}(m)$. 

Each $x\in\g$ acts on $\C^{a_j}$ by multiplication by $i\th_j$ for some 
$\th_j\in\R$, and is zero on $V$. Now $G$ is connected and abelian, so 
$G\cong\R^k$, and $\exp:\g\ra G$ is surjective. Hence we can write each 
$\ga\in G$ as $\exp(x)$ for $x\in\g$, and so $\ga$ acts on $\C^{a_j}$ by 
multiplication by ${\rm e}^{i\th_j}$, and as the identity on $V$. Putting 
$a_{n+1}=0$ and $\C^{a_{n+2}}=V$, we have shown that \eq{sy9eq1} satisfies 
(i), (v) and the second part of~(iii).

Now consider the general case with $G\subset\SU(m)\lt\C^m$. By
projecting $G$ from $\SU(m)\lt\C^m$ to $\SU(m)$ we can reduce to
the previous case, and decompose $\C^m$ into eigenspaces. However,
we now have to allow for $\g$ to act by {\it translations} in each 
factor, as well as by $\su(m)$ rotations.

Since the $\su(m)$ part of $\g$ acts on $\C^{a_1},\ldots,\C^{a_n}$ 
with nonzero eigenvalue, and $\g$ is abelian, by moving the origin 
in $\C^{a_j}$ we can eliminate the translation part, so that $G$ 
acts on $\C^{a_j}$ by multiplication by ${\rm e}^{i\th_j}$ for 
$j=1,\ldots,n$, as in (iii). Moving the origin is allowed, as we
seek only an {\it affine} isomorphism $\C^m\cong\C^{a_1}\t\cdots
\t\C^{a_{n+2}}$, rather than a vector space isomorphism.

However, moving the origin in $V$ in \eq{sy9eq1} has no effect 
on the $V$ translation-component of the action of $\g$, because 
this is the zero eigenspace of the $\su(m)$ part of $\g$. So we
cannot eliminate translations in the $V$ directions by choosing
the origin appropriately. Instead, define $\C^{a_{n+1}}$ to be the
complex vector subspace of $V$ generated by the $V$ 
translation-components of $\g$, and let $\C^{a_{n+2}}$ be the
orthogonal complement to $\C^{a_{n+1}}$ in~$V$.

Then we have an affine isomorphism $\C^m\cong\C^{a_1}\t\cdots
\t\C^{a_{n+2}}$ such that each $\ga\in G$ acts on $\C^{a_j}$ by 
multiplication by ${\rm e}^{i\th_j}$ for $j=1,\ldots,n$, $G$ 
acts by translations on $\C^{a_{n+1}}$, and $G$ acts trivially on 
$\C^{a_{n+2}}$, so that (i), (v) and parts of (iii) and (iv) are
satisfied.

Now we have put the $G$-action in a standard form, we prove the
`only if' part of the theorem. Suppose $L$ is special Lagrangian 
in $\C^m$, and $\phi(x)$ is normal to $L$ at $L$ for every $x$ in 
$\g$. The key idea we shall use is that for each $x$ in $\g$, as 
$\phi(x)$ is normal to $L$ and $L$ is Lagrangian, the vector field 
$I\bigl(\phi(x)\bigr)$ is tangent to $L$ at $L$. By exponentiating 
$I\bigl(\phi(x)\bigr)$ we get a 1-parameter family of diffeomorphisms 
of $\C^m$, which {\it locally} preserve $L$. That is, for each $z$ in 
the interior of $L$, there exists $\ep>0$ such that 
$\exp\bigl(t\,I(\phi(x))\bigr)z\in L$ for~$t\in(-\ep,\ep)$.

As $\g$ is abelian and the vector fields $\phi(\g)$ are holomorphic,
we see that
\begin{equation*}
\exp\bigl(I(\phi(x))\bigr)\circ\exp\bigl(I(\phi(y))\bigr)=
\exp\bigl(I(\phi(x+y))\bigr)
\end{equation*}
for all $x,y\in\g$. Thus the $\exp\bigl(I(\phi(x))\bigr)$ form
an abelian Lie group $\exp\bigl(I(\phi(\g))\bigr)$ of 
diffeomorphisms of $\C^m$, isomorphic to $\R^k$. We use this to 
extend $L$ to a {\it globally} invariant submanifold $L'$. Define
\begin{equation*}
L'=\bigcup_{x\in{\mathfrak g}}\exp\bigl(I(\phi(x))\bigr)(L).
\end{equation*}
Then it is not difficult to show that $L'$ is a special 
Lagrangian submanifold of $\C^m$ containing $L$, invariant 
under~$\exp\bigl(I(\phi(\g))\bigr)$.

One way to prove this is to use real analyticity, and the results 
of \S\ref{sy3}. Each connected component $L_i'$ of the interior 
of $L'$ contains a connected component $L_i$ of the interior of 
$L$. As $L_i$ is real analytic and $L_i'$ is the orbit of $L_i$ 
under a Lie group, $L_i'$ is also real analytic. So as the special
Lagrangian condition holds on a nonempty open subset $L_i$ of 
$L_i'$, it holds on all of~$L_i'$. 

As $L'$ is Lagrangian we have $\om\vert_{L'}\equiv 0$. 
But $\exp\bigl(I(\phi(x))\bigr)(L')=L'$ for $x\in\g$. Hence 
\e
\exp\bigl(I(\phi(x))\bigr)^*(\om)\big\vert_{L'}\equiv 0
\qquad\text{for all $x\in\g$.}
\label{sy9eq2}
\e
Write $\om=\sum_{j=1}^{n+2}\om^j$, where $\om^j$ is the
projection of $\om$ to $\C^{a_j}$. Let $x\in\g$ act on $\C^{a_j}$ 
by multiplication by $i\th_j$ for $j=1,\ldots,n$. Then
\e
\exp\bigl(I(\phi(x))\bigr)^*(\om)=
{\rm e}^{-2\th_1}\om^1+\cdots+
{\rm e}^{-2\th_n}\om^n+\om^{n+1}+\om^{n+2}.
\label{sy9eq3}
\e

Combining \eq{sy9eq2} and \eq{sy9eq3} for all $x\in\g$, and
remembering that the eigenvalues of $\g$ on $\C^{a_1},\ldots,\C^{a_n}$
are distinct and nonzero, we see that
\begin{equation*}
\om^j\vert_{L'}\equiv 0
\quad\text{for $j=1,\ldots,n$, and}\quad
\om^{n+1}\vert_{L'}+\om^{n+2}\vert_{L'}\equiv 0.
\end{equation*}
By considering the tangent spaces of $L'$ we find that $L'$
admits a local product structure, and deduce that 
$L'\subseteq L_1\t\cdots\t L_n\t N$, where $L_j$ is a
Lagrangian submanifold of $\C^{a_j}$ for $j=1,\ldots,n$, 
and $N$ is Lagrangian in $\C^{a_{n+1}}\t\C^{a_{n+2}}$. 

Let the $L_j$ and $N$ be as small as possible such that 
$L'\subseteq L_1\t\cdots\t L_n\t N$. This defines the
$L_j$ and $N$ uniquely. As $L'$ is special Lagrangian, 
it follows that $L_1,\ldots,L_n$ and $N$ are actually 
special Lagrangian in $\C^{a_1},\ldots,\C^{a_n}$ and 
$\C^{a_{n+1}}\t\C^{a_{n+2}}$ with some phases; we can 
fix the phases to be 1 by choosing the holomorphic volume 
forms on $\C^{a_j}$ appropriately. 

Since $L_1,\ldots,L_n$ and $N$ are defined uniquely using $L'$, 
which is invariant under $\exp\bigl(I(\phi(\g))\bigr)$, the 
$L_1,\ldots,L_n$ and $N$ must also be invariant under 
$\exp\bigl(I(\phi(\g))\bigr)$. But $\exp\bigl(I(\phi(x))\bigr)$ 
multiplies by $e^{-\th_j}$ in $\C^{a_j}$ for some $\th_j\in\R$ for 
$j=1,\ldots,n$, and $\th_j$ can take any value in $\R$ as $x$ varies. 
Therefore $L_j$ is invariant under all dilations of $\C^{a_j}$ for 
$j=1,\ldots,n$, and is a {\it cone} by Definition \ref{conedef}. 
This proves part (iii) of the theorem.

Similarly, $N$ is invariant under the action of 
$\exp\bigl(I(\phi(\g))\bigr)$ on $\C^{a_{n+1}}\t\C^{a_{n+2}}$.
Let $\cal O$ be an orbit of a point in $N$ under 
$\exp\bigl(I(\phi(\g))\bigr)$. Now $\exp\bigl(I(\phi(\g))\bigr)$ 
acts by translations on $\C^{a_{n+1}}$, and trivially on $\C^{a_{n+2}}$. 
Thus ${\cal O}=\R^l\t\{z\}$, where $\R^l$ is an affine subspace of
$\C^{a_{n+1}}$ and~$z\in\C^{a_{n+2}}$. 

By definition these translations  on $\C^{a_{n+1}}$ generate
$\C^{a_{n+1}}$ over $\C$, which forces $l\ge a_{n+1}$. However,
${\cal O}\subseteq N$ and $N$ is Lagrangian, so that 
$\om\vert_{\cal O}\equiv 0$. Therefore $l=a_{n+1}$, and 
$\R^l$ is a Lagrangian plane in $\C^{a_{n+1}}$. By choosing
the phase of the holomorphic volume form on $\C^{a_{n+1}}$
appropriately, we can assume that $\R^l$ is an SL
plane.

Thus $N$ is fibred by orbits $\cal O$ of the form $\R^{a_{n+1}}\t\{z\}$,
where $\R^{a_{n+1}}$ is an SL plane in $\C^{a_{n+1}}$
and $z\in\C^{a_{n+2}}$. It easily follows that $N=L_{n+1}\t L_{n+2}$,
where $L_{n+1}=\R^{a_{n+1}}$ in $\C^{a_{n+1}}$, and $L_{n+2}$ is 
an SL submanifold of $\C^{a_{n+2}}$. As 
$\exp\bigl(I(\phi(\g))\bigr)$ acts on $\C^{a_{n+1}}$ by 
translations in the direction of $\R^{a_{n+1}}$, it follows 
that $G=\exp\bigl(\phi(\g)\bigr)$ acts on $\C^{a_{n+1}}$ by 
translations in the direction of $I(\R^{a_{n+1}})$. 
This proves the `only if' part of Theorem \ref{sy9thm2}. But the 
`if' part follows very easily, given the discussion above, so the 
proof is complete.
\end{proof}

The theorem tells us that to apply Theorem \ref{sy9thm1}, we
need examples of SL cones $L_j$ in $\C^{a_j}$.
Now the most obvious SL cone in $\C^{a_j}$ is
$\R^{a_j}$. If we take $L_j=\R^{a_j}$ for all $j$ then $L'$ is
just $\R^m$ in $\C^m$, and we easily prove:

\begin{prop} Let\/ $2\le n\le m$, and let\/ $U(1)^{n-1}\t\R^{m-n}$ 
act on $\C^m$ by
\begin{align*}
({\rm e}^{i\th_1},&\ldots,{\rm e}^{i\th_{n-1}},x_{n+1},\ldots,x_m):
(z_1,\ldots,z_m)\longmapsto\\
&\bigl({\rm e}^{i\th_1}z_1,\ldots,{\rm e}^{i\th_{n-1}}z_{n-1},
{\rm e}^{-i(\th_1+\cdots+\th_{n-1})}z_n,z_{n+1}+ix_{n+1},\ldots,
z_m+ix_m\bigr),
\end{align*}
for $\th_1,\ldots,\th_{n-1}\in[0,2\pi)$ and\/ $x_{n+1},\ldots,x_m\in\R$.
Let\/ $G$ be any connected Lie subgroup of\/ $U(1)^{n-1}\t\R^{m-n}$, and\/ 
$L$ be $\R^m$ in $\C^m$. Then Theorem \ref{sy9thm1} applies to $G$ and\/ $L$, 
and constructs a family of\/ $G$-invariant SL submanifolds 
$N_c$ in $\C^m$ with phase $1$ or $i$, depending on~$c\in\g^*$.
\label{sy9prop}
\end{prop}

This gives many families of SL submanifolds in 
$\C^m$, which can be written down very explicitly. Here is an
example with~$G=\U(1)$.

\begin{ex} Let $a_1,\ldots,a_m$ be integers with $\hcf(a_1,\ldots,a_m)=1$
and $a_1+\cdots+a_m=0$, and let $G$ be $\U(1)$ acting on $\C^m$ by
\begin{equation*}
{\rm e}^{i\th}:(z_1,\ldots,z_m)\mapsto
\bigl({\rm e}^{ia_1\th}z_1,\ldots,{\rm e}^{ia_m\th}z_m\bigr).
\end{equation*}
Then $G$ lies in $\SU(m)$. The moment map of this $G$-action is
\begin{equation*}
\mu:(z_1,\ldots,z_m)\mapsto a_1\ms{z_1}+\cdots+a_m\ms{z_m}.
\end{equation*}
Take $L$ to be $\R^m$ in $\C^m$, and apply Proposition \ref{sy9prop}
and Theorem \ref{sy9thm1}. We find that for each $c\in\R$, the
subset $N^{a_1,\ldots,a_m}_c$ in $\C^m$ given by
\begin{align*}
N^{a_1,\ldots,a_m}_c=\Bigl\{
\bigl({\rm e}^{ia_1\th}x_1&,\ldots,{\rm e}^{ia_m\th}x_m\bigr):
\th\in[0,2\pi),\\ 
&x_1,\ldots,x_m\in\R,\qquad 
a_1x_1^2+\cdots+a_mx_m^2=c\Bigr\}
\end{align*}
is an SL $m$-fold with phase $i$. Define 
\begin{equation*}
{\cal H}^{a_1,\ldots,a_m}_c=\bigl\{(x_1,\ldots,x_m):x_j\in\R,\quad
a_1x_1^2+\cdots+a_mx_m^2=c\bigr\}.
\end{equation*}
Then $N^{a_1,\ldots,a_m}_c$ is the image of ${\cal S}^1\t
{\cal H}^{a_1,\ldots,a_m}_c$ under the map
\begin{equation*}
\Phi:\bigl({\rm e}^{i\th},(x_1,\ldots,x_m)\bigr)\mapsto
\bigl({\rm e}^{ia_1\th}x_1,\ldots,{\rm e}^{ia_m\th}x_m\bigr).
\end{equation*}
When $c\ne 0$ this is an {\it immersion}, so that 
$N^{a_1,\ldots,a_m}_c$ is a nonsingular immersed $m$-submanifold. 
Also $N^{a_1,\ldots,a_m}_0$ is a {\it cone} with an isolated
singular point at 0, and is otherwise nonsingular as an immersed
submanifold.

However, $\Phi$ is generally not injective. Generically $\Phi$ is 
2:1, since
\e
\Phi\bigl({\rm e}^{i\th},(x_1,\ldots,x_m)\bigr)=
\Phi\bigl(-{\rm e}^{i\th},((-1)^{a_1}x_1,\ldots,(-1)^{a_m}x_m)\bigr).
\label{sy9eq4}
\e
Thus we may regard $N^{a_1,\ldots,a_m}_c$ as an immersion of 
$({\cal S}^1\t{\cal H}^{a_1,\ldots,a_m}_c)/\Z_2$ in $\C^m$.
When two or more of the $x_j$ vanish $\Phi$ may become $2k:1$
for $k>1$, and then $N^{a_1,\ldots,a_m}_c$ is singular as an
embedded submanifold.
\label{cmex3}
\end{ex}

Here is some more detail on the topology of this example when~$m=3$.

\begin{ex} Let $a_1,a_2$ be positive, coprime integers and
$a_3=-a_1-a_2$, let $c\in\R$, and let $N^{a_1,a_2,a_3}_c$ 
and other notation be as in Example \ref{cmex3}. Then if 
$c>0$ then $\Phi$ is 2:1 everywhere, and $N^{a_1,a_2,a_3}_c$ 
is an embedded 3-fold diffeomorphic to $(T^2\t\R)/\Z_2$, where 
$\Z_2$ acts freely on $T^2\t\R$. If $a_3$ is even then $\Z_2$ acts 
trivially on $\R$, and $N^{a_1,a_2,a_3}_c$ is diffeomorphic
to $T^2\t\R$. If $a_3$ is odd, then $N^{a_1,a_2,a_3}_c$ is
the total space of a nontrivial real line bundle over the
Klein bottle, and has only one end modelled on~$T^2\t(0,\iy)$.

When $c=0$, there are two cases: if $a_3$ is even then 
$N^{a_1,a_2,a_3}_0$ is the union of two opposite, embedded 
$T^2$-cones, meeting at 0, and if $a_3$ is odd then 
$N^{a_1,a_2,a_3}_0$ is just one embedded $T^2$-cone. The 
difference is that ${\cal H}^{a_1,a_2,a_3}_0\sm\{0\}$ separates 
into two cones with $x_3>0$ and $x_3<0$, and whether $a_3$ is even
or odd determines whether the identity \eq{sy9eq4} fixes these
cones or swaps them.

When $c<0$, if $a_3$ is even then $N^{a_1,a_2,a_3}_c$ is the 
union of two immersed copies of ${\cal S}^1\t\R^2$, and if 
$a_3$ is odd then $N^{a_1,a_2,a_3}_c$ is just one immersed
${\cal S}^1\t\R^2$. Note also that $\Phi$ maps $-2a_3$ points 
of the form $\bigl({\rm e}^{i\th},(0,0,x_3)\bigr)$ to one point 
in $\C^3$, so $N^{a_1,a_2,a_3}_c$ is singular as an embedded 
submanifold along an ${\cal S}^1$ in~$\C^3$. 

It can be shown that $N^{1,1,-2}_c$ is isomorphic to $L_{c/2,c/2,0}$ 
of Example \ref{c3ex1} under a change of coordinates, and has a $T^2$ 
symmetry group. But in general $N^{a_1,a_2,a_3}_c$ is a new example, 
not of cohomogeneity one for~$c\ne 0$.
\label{c3ex4}
\end{ex}

Here is another example in $\C^3$.

\begin{ex} Let $G$ be $\R$, acting on $\C^3$ by
\begin{equation*}
t:(z_1,z_2,z_3)\mapsto
\bigl({\rm e}^{it}z_1,{\rm e}^{-it}z_2,z_3+it\bigr).
\end{equation*}
Then $G$ lies in $\SU(3)$. The moment map of this $G$-action is
\begin{equation*}
\mu:(z_1,z_2,z_3)\mapsto \ms{z_1}-\ms{z_2}+2\Re z_3.
\end{equation*}
Applying Proposition \ref{sy9prop} and Theorem \ref{sy9thm1}
with $c=0$ shows that
\begin{equation*}
N=\Bigl\{\bigl({\rm e}^{it}x_1,{\rm e}^{-it}x_2,x_3+it\bigr):
t,x_1,x_2,x_3\in\R,\quad x_1^2-x_2^2+2x_3=0\Bigr\}
\end{equation*}
is an SL 3-fold in $\C^3$ with phase $i$, which
is nonsingular and diffeomorphic to $\R^3$. One can picture $N$
as being a bit like a helicoid in~$\R^3$.
\label{c3ex5}
\end{ex}

In our last example we construct a family of SL 4-folds $N_c$ 
in $\C^4$ out of an SL cone $L_1$ in $\C^3$, using the ideas of 
Theorems \ref{sy9thm1} and~\ref{sy9thm2}.

\begin{ex} Let $G$ be $\U(1)$ acting on $\C^4$ by
\begin{equation*}
{\rm e}^{i\th}:(z_1,\ldots,z_4)\mapsto
\bigl({\rm e}^{i\th}z_1,{\rm e}^{i\th}z_2,{\rm e}^{i\th}z_3,
{\rm e}^{-3i\th}z_4\bigr).
\end{equation*}
Then $G\subset\SU(4)$, and has moment map
\begin{equation*}
\mu:(z_1,\ldots,z_4)\mapsto \ms{z_1}+\ms{z_2}+\ms{z_3}-3\ms{z_4}.
\end{equation*}
Let $L_1$ be an SL cone in $\C^3$. Then applying
Theorems \ref{sy9thm1} and \ref{sy9thm2} to $G$ and $L=L_1\t\R$ 
in $\C^3\t\C$, we find that for each $c\in\R$,
\begin{align*}
N_c=\Bigl\{\bigl({\rm e}^{i\th}x_1&,{\rm e}^{i\th}x_2,
{\rm e}^{i\th}x_3,{\rm e}^{-3i\th}x_4\bigr):
\th\in[0,2\pi),\\ 
&(x_1,x_2,x_3)\in L_1,\quad 
x_4\in\R,\quad
\ms{x_1}+\ms{x_2}+\ms{x_3}-3x_4^2=c\Bigr\}
\end{align*}
is an SL 4-fold in $\C^4$ with phase~$i$.

Now let $L_1$ be a cone on a compact Riemann surface $\Si$, with 
an isolated singular point at 0, and suppose for simplicity that
$L_1\cap{\rm e}^{2\pi i/3}L_1=\{0\}$. Then for $c>0$ we find
that $N_c$ is a nonsingular, embedded SL 4-fold 
diffeomorphic to $\Si\t{\cal S}^1\t\R$. Similarly, $N_0$ is a cone 
on 2 $\Si\t{\cal S}^1$ with an isolated singularity at 0, and if 
$c<0$ then $N_c$ is singular on an ${\cal S}^1$ in~$\C^4$.
\label{c4ex2}
\end{ex}

\end{document}